%
\def\date{19 December  2018}  
\magnification=1200
\overfullrule=0pt
\input epsf.tex
\input miniltx %
\input graphicx.sty %
\input color
\newif\ifproofmode
\def\xrefsfilename{ext.xrf}  
\def\myinput#1{\immediate\openin0=#1\relax
   \ifeof0\write16{Cannot input file #1.}
   \else\closein0\input#1\fi}
\newcount\referno
\newcount\thmno
\newcount\secno
\newcount\figno
\referno=0
\thmno=0
\secno=0
\def\ifundefined#1{\expandafter\ifx\csname#1\endcsname\relax}
\myinput \xrefsfilename
\immediate\openout1=\xrefsfilename
\def\bibitem#1#2\par{\ifundefined{REFLABEL#1}\relax\else
 \global\advance\referno by 1\relax
 \immediate\write1{\noexpand\expandafter\noexpand\def
 \noexpand\csname REFLABEL#1\endcsname{\the\referno}}
 \global\expandafter\edef\csname REFLABEL#1\endcsname{\the\referno}
 \item{\the\referno.}#2\ifproofmode [#1]\fi\fi}
\def\cite#1{\ifundefined{REFLABEL#1}\ignorespaces
   \global\expandafter\edef\csname REFLABEL#1\endcsname{?}\ignorespaces
   \write16{ ***Undefined reference #1*** }\fi
 \csname REFLABEL#1\endcsname}
\def\nocite#1{\ifundefined{REFLABEL#1}\ignorespaces
   \global\expandafter\edef\csname REFLABEL#1\endcsname{?}\ignorespaces
   \write16{ ***Undefined reference #1*** }\fi}
\def\newthm#1#2\par{\global\advance\thmno by 1\relax
 \immediate\write1{\noexpand\expandafter\noexpand\def
 \noexpand\csname THMLABEL#1\endcsname{(\the\secno.\the\thmno)}}
 \global\expandafter\edef\csname THMLABEL#1\endcsname{(\the\secno.\the\thmno)}
 \bigbreak\penalty-500\noindent{\bf(\the\secno.\the\thmno)\enspace}\ignorespaces
 \ifproofmode {\bf[#1]} \fi{\sl#2}
 \medbreak\penalty-200}
\def\newsection#1#2\par{\global\advance\secno by 1\relax
 \immediate\write1{\noexpand\expandafter\noexpand\def
 \noexpand\csname SECLABEL#1\endcsname{\the\secno}}
 \global\expandafter\edef\csname SECLABEL#1\endcsname{\the\secno}
 \vskip0pt plus.3\vsize
 \vskip0pt plus-.3\vsize\bigskip\bigskip\vskip\parskip\penalty-250
 \message{\the\secno. #2}\thmno=0
 \centerline{\bf\the\secno. #2\ifproofmode {\rm[#1]} \fi}
 \nobreak\smallskip\noindent}
\def\refthm#1{\ifundefined{THMLABEL#1}\ignorespaces
 \global\expandafter\edef\csname THMLABEL#1\endcsname{(?)}\ignorespaces
 \write16{ ***Undefined theorem label #1*** }\fi
 \csname THMLABEL#1\endcsname}
\def\refsec#1{\ifundefined{SECLABEL#1}\ignorespaces
 \global\expandafter\edef\csname SECLABEL#1\endcsname{(?)}\ignorespaces
 \write16{ ***Undefined section label #1*** }\fi
 \csname SECLABEL#1\endcsname}
\def\newfig#1{\global\advance\figno by 1\relax
 \immediate\write1{\noexpand\expandafter\noexpand\edef
 \noexpand\csname FIGLABEL#1\endcsname{\the\figno}}\ignorespaces
 \global\expandafter\edef\csname FIGLABEL#1\endcsname{\the\figno}\ignorespaces
 \the\figno\ifproofmode{\bf[#1]}\fi}
\def\reffig#1{\ifundefined{FIGLABEL#1}\ignorespaces
 \global\expandafter\edef\csname FIGLABEL#1\endcsname{??}\ignorespaces
 \write16{ ***Undefined figure label #1*** }\fi
 \csname FIGLABEL#1\endcsname}
\def\xx#1{{\color{red}#1}}
\def\ifc{internally $4$-connected}
\def\afc{almost $4$-connected}
\def\afcity{almost $4$-connectivity}
\def\ifcity{internal $4$-connectivity}
\def\he{homeomorphic embedding}
\def\emb{\hookrightarrow}

\def\calc{{\cal C}}

\def\stable{rigid}
\font\smallrm=cmr8
\def\junk#1{}
\def\dfn#1{{\sl#1}}
\def\cond#1#2\par{\smallbreak\noindent\rlap{\rm(#1)}\ignorespaces
\hangindent=30pt\hskip30pt{\rm#2}\smallskip}
\def\claim#1#2\par{{\medbreak\noindent\rlap{\rm(#1)}\ignorespaces
  \rightskip20pt\hangindent=20pt\hskip20pt{\ignorespaces\sl#2}\smallskip}}
\def\leanclaim#1#2\par{{\medbreak\noindent\rlap{\rm(#1)}\ignorespaces
  \hangindent=20pt\hskip20pt{\ignorespaces\sl#2}\smallskip}}
\def\proof{\smallbreak\noindent{\sl Proof. }}

\def\qed{\hfill$\square$\bigskip\medskip}
\def\sqr#1#2{{\vcenter{\vbox{\hrule height.#2pt
\hbox{\vrule width.#2pt height #1pt \kern#1pt
\vrule width.#2pt}
\hrule height.#2pt}}}}
\def\square{\mathchoice\sqr56\sqr56\sqr{2.1}3\sqr{1.5}3}
\outer\def\beginsection#1\par{\vskip0pt plus.3\vsize
   \vskip0pt plus-.3\vsize\bigskip\bigskip\vskip\parskip
   \message{#1}\centerline{\bf#1}\nobreak\smallskip\noindent}
\newcount\remarkno
\def\REMARK#1{{%
   \footnote{${}^{\the\remarkno}$}{\baselineskip=11pt #1
   \vskip-\baselineskip}\global\advance\remarkno by1}}

\nocite{KelVertex}
\nocite{KelSpecialK4}
\nocite{Kelcfc}
\nocite{KelPlanarity}
\nopagenumbers
\footline={\hfil}
\baselineskip=12pt
\phantom{a}\vskip .25in
\centerline{{\bf NON-PLANAR EXTENSIONS OF SUBDIVISIONS OF PLANAR GRAPHS}}
\vskip.5in
\bigskip
\centerline{{\bf Sergey Norin}$^{1}$\vfootnote{$^1$}
{\smallrm Supported by an NSERC discovery grant.}
}
\centerline{Department of Mathematics and Statistics}
\centerline{McGill University}
\centerline{Montreal, Quebec H3A 2K6, Canada}
\bigskip
\centerline{and}
\bigskip
\centerline{{\bf Robin Thomas}%
$^{2}$\vfootnote{$^2$}{\smallrm Partially
supported
by NSF under Grants No.~DMS-9623031, DMS-0200595 and DMS-1202640, and by NSA under
Grant No.~MDA904-98-1-0517.
}}
\centerline{School of Mathematics}
\centerline{Georgia Institute of Technology}
\centerline{Atlanta, Georgia  30332-0160, USA}


\vfill
\baselineskip 11pt
\noindent 5 November 1998,
Revised \xx{\date}.
\hfil\break\noindent 
\xx{Published in {\it  J.~Combin.\ Theory Ser.~B \bf121} (2016), 326-366.
This version fixes an error in the published paper. The error was kindly pointed out to us by Katherine Naismith.
Changes from the published version are indicated in red.
}%
\vfil\eject
\baselineskip 18pt
\footline{\hss\tenrm\folio\hss}

\beginsection ABSTRACT

\dfn{Almost $4$-connectivity} is a weakening of $4$-connectivity which allows for vertices of degree three.
In this paper we prove the following theorem. Let $G$ be an \afc\ triangle-free planar graph,
and let $H$ be an \afc\
non-planar graph such that $H$ has a subgraph isomorphic to a subdivision
of $G$. Then there exists a graph $G'$ such that $G'$ is isomorphic
to a minor of $H$, and either
\item{(i)}$G'=G+uv$ for some vertices $u,v\in V(G)$ such that
no facial cycle of $G$ contains both $u$ and $v$, or
\item{(ii)}$G'=G+u_1v_1+u_2v_2$ for some distinct vertices
$u_1,u_2,v_1,v_2\in V(G)$ such that $u_1,u_2,v_1,v_2$ appear on
some facial cycle of $G$ in the order listed.

\noindent
This is a lemma to be used in other papers. In fact, we prove a more general theorem, where we relax the connectivity
assumptions, do not assume that $G$ is planar, and consider subdivisions rather than minors. Instead of face boundaries we work with a collection of cycles that cover every edge twice and have pairwise connected intersection.
Finally, we prove a version of this result that applies when $G\backslash X$ is planar for some set $X\subseteq V(G)$ of size at most $k$, but $H\backslash Y$ is non-planar for every set $Y\subseteq V(H)$ of size at most $k$.

\vfil\eject

\newsection{intro}INTRODUCTION

In this paper graphs are finite and simple (i.e., they have no
loops or multiple edges).
\dfn{Paths} and \dfn{cycles} have no ``repeated" vertices or edges.
A graph is a \dfn{subdivision} of another
if the first
can be obtained from the second by replacing each edge by a non-zero length
path
with the same ends, where the paths are disjoint, except possibly
for shared ends.  The replacement paths are called \dfn{segments}, and
their ends are called \dfn{branch-vertices}.
For later convenience a one-vertex component of a graph is
also regarded as a segment, and its unique vertex as a branch-vertex.
Let $G,S,H$ be graphs such that $S$ is a subgraph of $H$ and is isomorphic to a
subdivision of $G$.  In that case we say that $S$ is a $G$-\dfn{subdivision}
in $H$.  If $G$ has no vertices of degree two (which will be the case in
our applications), then the segments and branch-vertices of $S$ are
uniquely determined by $S$.
An \dfn{$S$-path} is a path of length at least one with both ends in
$S$ and otherwise disjoint from $S$.
A graph $G$ is \dfn{\afc} if it is simple, $3$-connected, has at least five
vertices, and $V(G)$ cannot be partitioned into three sets
$A,B,C$ in such a way that $|C|=3$, $|A|\ge2$, $|B|\ge2$, and no
edge of $G$ has one end in $A$ and the other end in $B$.

Let a non-planar graph $H$ have a subgraph $S$ isomorphic to a subdivision
of a planar graph $G$. For various problems in structural graph
theory it is useful to
know  the minimal subgraphs of $H$
that have a subgraph isomorphic to a subdivision of $G$ and are non-planar.
We show that under some mild connectivity assumptions these
``minimal non-planar extensions" of $G$ are quite nice:

\newthm{main}
Let $G$ be an \afc\ planar graph on at least seven vertices,
let $H$ be an \afc\
non-planar graph, and let there exist a $G$-subdivision in $H$.
Then there exists a $G$-subdivision
$S$ in $H$ such that one of the following
conditions holds:
\item{(i)}there exists an $S$-path in $H$
joining two vertices of $S$ not incident with the same face, or
\item{(ii)}there exist two disjoint $S$-paths with ends
$s_1,t_1$ and $s_2,t_2$, respectively, such that the vertices
$s_1,s_2,t_1,t_2$ belong to some face boundary of $S$ in
the order listed. Moreover, for $i=1,2$ the vertices $s_i$ and $t_i$ do not
belong to the same segment of $S$, and if two segments of $S$ include all
of $s_1,t_1,s_2,t_2$, then those segments are vertex-disjoint.


\noindent The connectivity assumptions guarantee that the face
boundaries in a planar embedding of $S$ are uniquely determined, and
hence it makes sense to speak about incidence with faces.
Theorem~\refthm{main}  is related to, but independent of [\cite{RobSeyThoCubic}].
We refer the reader to [\cite{ThoSurvey}] for an overview of related results.

In Section~\refsec{planar} we deduce the following corollary, stated
there as~\refthm{minor}.
A graph is a \dfn{minor} of another if the first can be obtained from
a subgraph of the second by contracting edges.
If $G$ is a graph and $u,v\in V(G)$ are not adjacent, then
by $G+uv$ we denote the graph obtained from $G$ by adding an edge
with ends $u$ and $v$.

\newthm{maincor}Let $G$ be an \afc\ triangle-free planar graph,
and let $H$ be an \afc\
non-planar graph such that $H$ has a subgraph isomorphic to a subdivision
of $G$. Then there exists a graph $G'$ such that $G'$ is isomorphic
to a minor of $H$, and either
\item{(i)}$G'=G+uv$ for some vertices $u,v\in V(G)$ such that
no facial cycle of $G$ contains both $u$ and $v$, or
\item{(ii)}$G'=G+u_1v_1+u_2v_2$ for some distinct vertices
$u_1,u_2,v_1,v_2\in V(G)$ such that $u_1,u_2,v_1,v_2$ appear on
some facial cycle of $G$ in the order listed.

\noindent
While the statement of~\refthm{maincor} is nicer, it has the drawback
that we assume that $H$ has a subgraph isomorphic to a {\sl subdivision} of $G$,
and deduce that it has only a {\sl minor} isomorphic to $G'$.
That raises the question whether there is a similar theory that applies
when $H$ has a minor isomorphic to $G$.
Such a theory indeed exists and is developed in~[\cite{HegTho}],
using~\refthm{rajneesh} below.
Informally, there is an analogue of~\refthm{main}, where either of the
two outcomes may be preceded by up to two vertex splits (inverse operations
to edge contraction).

In the applications of~\refthm{main} the graph $G$ is known explicitly,
but $H$ is not, and we are trying to deduce some information about
$H$. Since it is possible to generate all graphs that can be obtained
from subdivisions of $G$ by means of~\refthm{main}(i) or~\refthm{main}(ii),
we thus obtain a list of specific
non-planar graphs such  that $H$ has
a subgraph isomorphic to a subdivision of one of the graphs in the list.
The graphs $G$ of interest in applications tend to possess
a lot of symmetry, and so the generation process is usually less daunting
than it may seem.

A sample application of our result is presented in Section~\refsec{appl},
but let us informally describe the applications
from~[\cite{DinOpoThoVerlarge}]
and~[\cite{ThoThoTutte}].
Theorem~\refthm{polyplanar}, a close relative of~\refthm{main},
 is used in [\cite{DinOpoThoVerlarge}] to show that
for every positive integer $k$, there is an integer $N$
such that every $4$-connected non-planar graph with at least
$N$ vertices has a minor isomorphic to the complete bipartite graph $K_{4,k}$,
or the graph obtained from a cycle of length $2k+1$ by adding an edge
joining every pair of vertices at distance exactly $k$, or
the graph obtained from a cycle of length $k$ by adding two vertices
adjacent to
each other and to every vertex on the cycle.
Using this Bokal, Oporowski, Richter and Salazar~[\cite{BokOpoRicSal}]
proved that,
except for one well-defined infinite family, there
are only finitely many graphs of crossing number at least two
that are minimal in a specified sense.

In [\cite{ThoThoTutte}] it is shown that every \afc\ non-planar graph of girth
at least five has a subgraph isomorphic to
a subdivision of $P_{10}^-$,
the Petersen graph with one edge deleted. (It follows from
this that Tutte's $4$-flow conjecture [\cite{Tut66}] holds for graphs with
no subdivision isomorphic to $P_{10}^-$.) The way this is done is that
first it is shown that if $G$ is a graph of girth at least five and
minimum degree at least three, then it has a subgraph isomorphic to
a subdivision of
the Dodecahedron or $P_{10}^-$. Corollary~\refthm{maincor} is then used
to show that if $G$ is an \afc\ non-planar graph with a subgraph isomorphic
to a subdivision of the Dodecahedron, then $G$ has a subgraph isomorphic to
a subdivision of $P_{10}^-$.

We actually prove several results that are more general than~\refthm{main}.
It turns out that global planarity is not needed for the proof to go
through; thus we formulate most of our results in terms of not necessarily
planar graphs with a specified set $\cal C$ of cycles that cover every
edge twice, have pairwise empty or connected intersection, and satisfy another
natural condition.
We call such sets of cycles disk systems.
To deduce~\refthm{main} we let $\cal C$ be the disk system of facial
cycles in $S$.
This greater generality allows us to prove an analogue of~\refthm{main}
for graphs on higher surfaces.

We also investigate an extension of our original problem to apex graphs.
What can we say when $G$ has a set $X\subseteq V(G)$ of size at most
$k$ such that $G\backslash X$ is planar, but $H$ has no such set?
Is there still an analogue of~\refthm{main}?
To prove an exact analogue seems to be a difficult problem that will
require a complicated answer.
Luckily, for our applications we can assume that $G$ is triangle-free,
and we can afford to ``sacrifice" a few edges from $X$ to $G\backslash X$.
With those two simplifying assumptions we were able to prove~\refthm{apexcor},
a result along the lines of~\refthm{main}, that is simple enough to allow
a concise statement and yet strong enough to allow us to deduce the
desired applications.
One such application can be found in~[\cite{KawNorThoWolbdtw}].

The paper is organized as follows.
Throughout the paper we will have to transform one $G$-subdivision to another,
and it will be useful to keep track of the changes we have (or have not)
made.
There are four kinds of such transformations, called reroutings,
and we introduce them in Section~\refsec{reroute}.
In Section~\refsec{stable} we prove
a useful and well-known lemma which says that if a graph $H$ has a subgraph
isomorphic to a subdivision of a graph $G$ and $H$ is
$3$-connected, then $H$ has a subgraph isomorphic to a
subdivision of $G$ such that all ``bridges" are ``\stable".
In fact, we need a version of this for graphs that are not necessarily
$3$-connected.
We also review several basic results about planar graphs in
Section~\refsec{stable}.
In Section~\refsec{disksys} we introduce disk systems and
prove a version of our main result
without assuming any connectivity of $G$ or $H$.
In Section~\refsec{triads} we eliminate one of the outcomes by assuming
that $H$ is \afc, and in Section~\refsec{planar} we prove~\refthm{main}
and a couple of closely related theorems.
In Section~\refsec{appl} we
illustrate the use of~\refthm{main}.
Section~\refsec{improve} contains a technical improvement of one
of the earlier lemmas for use in Section~\refsec{apex}, where
we prove a version of our result when $G$ is at most $k$ vertices
away from being planar and $H$ is not. In Section~\refsec{pinwheel} we present an application
of this version of the result. 

\newsection{reroute} REROUTINGS

We will need a fair amount of different kinds of reroutings that transform
one $G$-subdivision into another, and in
order to avoid confusion it seems best to collect them all in one
place for easy reference.
If $P$ is a path and $x,y\in V(P)$, then $xPy$ denotes the subpath
of $P$ with ends $x$ and $y$.

First we recall the classical notion of a bridge.
Let $S$ be a subgraph of a graph $H$.  An $S$-\dfn{bridge} in $H$ is a
connected subgraph $B$ of $H$ such that $E(B)\cap E(S)=\emptyset$
and either $E(B)$ consists of a unique edge with both ends in $S$, or for
some component $C$ of $H\backslash V(S)$ the set $E(B)$ consists of all
edges of $H$ with at least one end in $V(C)$.  The vertices in
$V(B)\cap V(S)$ are called the \dfn{attachments} of $B$.

 Let $G,H$ be graphs, let $G$ have no vertices of degree two,
let $S$ be a $G$-subdivision in $H$,
let $v$ be a vertex of $S$ of degree $k$, let
$P_1,P_2,\ldots,P_k$ be the segments of $S$
incident with $v$, and let their other ends be $v_1,v_2,\ldots,v_k$,
respectively. Let $x,y\in V(P_1\cup P_2\cup\ldots P_k)$ be
distinct vertices, and let $Q$ be an $S$-path with ends $x$ and $y$.
Furthermore, let $P$ be a suitable subpath of $S$, to be specified later.
We wish to define a new $G$-subdivision $S'$ by removing all edges and
internal vertices of $P$ from $S\cup Q$.
 If $x,y\in V(P_1)$, $P_1$ has length at least two and
$P=xP_1y$, then we say that $S'$ is obtained from $S$ by
an \dfn{{\rm I}-rerouting}.
If, in addition, the $S$-bridge containing $Q$ has all attachments
in $P_1$, then we say that $S'$ is obtained from $S$ by
a \dfn{proper {\rm I}-rerouting}. See Figure~\reffig{ireroute}.
We emphasize that we indeed require that $P_1$ have at least two edges.

\bigskip

\bigskip
\centerline{\includegraphics[scale=1.50]{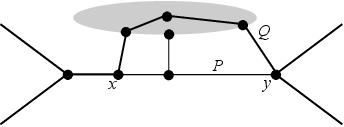}}
\bigskip
\centerline{Figure~\newfig{ireroute}. Proper I-rerouting.}
\bigskip

Let $k=3$, let $x\in V(P_1)-\{v\}$,  let $y$ be an internal
vertex of $P_2$,
and let $P=xP_1v$.
In those circumstances we say that $S'$ is obtained from $S$ by
a T-\dfn{rerouting}. See Figure~\reffig{treroute}.


\bigskip
\centerline{\includegraphics[scale=1.70]{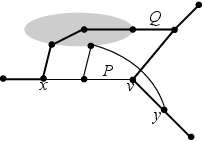}}
\bigskip
\centerline{Figure~\newfig{treroute}. T-rerouting.}
\bigskip

If $k\ge4$ and there exists an integer $i\in\{1,2\}$ such that
$P_i$ has length at least two,
$x,y\in V(P_i)$ and $P=xP_iy$, then we say that $S'$ is obtained from $S$
by a \dfn{{\rm V}-rerouting}, and we say that it is obtained by a
\dfn{proper {\rm V}-rerouting} if all the attachments of
the $S$-bridge containing $Q$ belong to $P_1\cup P_2$.
In that case we say that $S'$ is obtained from $S$ by a
\dfn{proper {\rm V}-rerouting based at} $P_1$ and  $P_2$.
Thus a V-rerouting is also an I-rerouting,
but not so for proper reroutings.


The last type of rerouting which we define in this paragraph differs from all the types defined so far, as we remove the interiors of two paths from $S$ rather than one.
 Let $k\ge4$, let $x_1,x_2\in V(P_1)-\{v\}$
and $y_1,y_2\in V(P_2)-\{v\}$ be distinct vertices such that the vertices
 $x_1,x_2,v,y_1,y_2$ appear on the
path $P_1\cup P_2$ in the order listed, and for $i=1,2$
let $Q_i$ be an $S$-path in $H$ with ends $x_i$ and $y_i$
such that $Q_1$ and $Q_2$ are disjoint.
 Let $S'$ be obtained from $S\cup Q_1\cup Q_2$ by deleting the edges
and internal vertices of the paths $x_1P_1x_2$ and $y_1P_2y_2$.
Then $S'$ is a
$G$-subdivision in $H$, and we say that $S'$ is obtained from $S$ by
an \dfn{{\rm X}-rerouting} of $S$.
See Figure~\reffig{xreroute}.
It is obtained by
a \dfn{proper {\rm X}-rerouting}
if the bridges containing $Q_1$ and $Q_2$ have all their attachments
in $P_1\cup P_2$.
In that case we say that $S'$ is obtained from $S$ by a
\dfn{proper {\rm X}-rerouting based at} $P_1$ and  $P_2$.
We say that $S'$ is obtained from $S$ by a \dfn{rerouting} if it is
obtained from $S$ by an I-rerouting,
a V-rerouting, a T-rerouting or an X-rerouting.
This relation is not symmetric, because in an I-rerouting and V-rerouting
we require that the path that is being changed have length at least two.
The distinction among different kinds of rerouting as well as proper reroutings
will not be needed until the last two sections and may be safely ignored
until then.

\bigskip

\bigskip
\centerline{\includegraphics[scale=2.0]{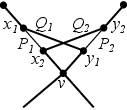}}
\bigskip
\centerline{Figure~\newfig{xreroute}. X-rerouting.}
\bigskip

\junk{
We have defined distant attachment of an I-rerouting above; for the other
reroutings we define distant attachment as follows.
We say that a vertex $w\in V(S)$ is a \dfn{distant attachment} of
the rerouting if  $w$ is an attachment of the $S$-bridge containing $Q$
(or one of the two $S$-bridges containing $Q_1$ or $Q_2$ in case of
an X-rerouting).
This creates a possible ambiguity, because an I-rerouting could also be
a V-rereouting or a T-rerouting, and the definition of distant attachment
depends on whether the rerouting is regarded as an I-rerouting or not.
To clarify this we say that the I-rerouting version of the definition will
only be applied when we explicitly refer to the rerouting as an I-rerouting.
} 

\newsection{stable}RIGID BRIDGES AND PLANAR GRAPH LEMMAS



Let $G$ be a graph with no vertices of degree two, and let
$S$ be a $G$-subdivision in a graph $H$.
If $B$ is an $S$-bridge of $H$, then we
say that $B$ is \dfn{unstable} if it has at least one attachment and
some segment of $S$ includes all the
attachments of $B$; otherwise we say that $B$ is
\dfn{rigid}.
Our next lemma, essentially due to Tutte,
says that in a $3$-connected graph it is possible to make all
bridges rigid by changing $S$ using proper I-rerouting only.
A \dfn{separation} of a graph $G$ is a pair
$(A,B)$ of subsets of $V(G)$ such that $A\cup B = V(G)$, and there is
no edge between $A-B$ and $B-A$.  The \dfn{order} of $(A,B)$
is $|A\cap B|$.
We say that an $S$-bridge $J$ is \dfn{$2$-separated} from $S$ if there
exists a segment $Z$ of $S$, two (not necessarily distinct) vertices
$u,v\in V(Z)$ and a separation $(A,B)$ of $H$ such that
$A$ includes all branch-vertices of $S$,
$V(J\cup uZv)\subseteq B$ and $A\cap B=\{u,v\}$.




\newthm{stable}Let $G$ be a graph with no
vertices of degree two, let $H$ be a graph,
and let $S$ be a $G$-subdivision in $H$.
Then there exists a $G$-subdivision $S'$ in $H$
 obtained from $S$ by a sequence of proper {\rm I}-reroutings
such that every unstable $S'$-bridge is $2$-separated from $S'$.

\proof
We use the same argument as in~[\cite{KawNorThoWolbdtw}, Lemma~2.1],
but we give the proof for completeness.
We may choose a $G$-subdivision $S'$ obtained from $S$ by a
sequence of proper I-reroutings such that the number of edges that
belong to \stable\ $S'$-bridges is maximum.  We will show that $S'$ is as
desired.  To that end we may assume that $S'$ has a segment $Z$ such
that some $S'$-bridge that has at least one attachment has
all its attachments in $Z$.

Let $v_0, v_1,\dots, v_k$ be distinct vertices of $Z$, listed in order of
occurrence on $Z$ such that $v_0$ and $v_k$ are the ends of $Z$ and
$\{v_1,\dots, v_{k-1}\}$ is the set of all internal vertices of $Z$
that are attachments of a rigid $S'$-bridge.
We may assume that $k\ge 2$, for otherwise every $S'$-bridge with all
attachments in $Z$ is $2$-separated from $S'$.
Now let $J$ be an
$S'$-bridge with at least one attachment and all
attachments contained in $Z$, and let
$x,y\in V(Z)$ be the two (not necessarily distinct)
attachments of $J$ that maximize  $xZy$.
We claim that for $i=1,2,\ldots,k-1$ the vertex $v_i$ does not belong
to the interior of $xZy$.
To prove this claim suppose to the contrary that $v_i$ belongs to the
interior of $xZy$.
Then replacing the path $xZy$ by a subpath of $J$ with  ends $x$ and $y$
is a proper I-rerouting that produces a $G$-subdivision $S''$
with strictly more edges belonging to \stable\ $S''$-bridges,
because every edge that belongs to a \stable\ $S'$-bridge belongs
to a \stable\ $S''$-bridge, and both edges of $S'$ incident with $v_i$ belong
to a \stable\ $S''$-bridge,
contrary to the choice of $S'$.
This proves our claim that $v_i$ does not belong to the interior of $xZy$.
Thus there exists an integer $i=1,2,\ldots,k$ such that
$xZy$ is a subpath of $v_{i-1}Zv_i$.
Let $B$ be the union of the vertex-set of $v_{i-1}Zv_i$ and the
vertex-sets of all unstable $S'$-bridges whose attachments are contained in
$v_{i-1}Zv_i$,
and let $A:=V(H)-(B-\{v_{i-1},v_i\})$.
Then the earlier claim implies that $(A,B)$ is a separation,
witnessing that $J$ is $2$-separated from $S$, as desired.~\qed


We will need the following result, a relative
of  [\cite{KelSpecialK4}, \cite{RobSeyGM9},
\cite{SeyDisj}, \cite{Shi}, \cite{Tho2link}].
If $G$ is a graph and $X\subseteq V(G)$, then $G[X]$ denotes
the graph $G\backslash (V(G)-X)$.

\newthm{2paths}Let $G$ be a graph, and let $C$ be a cycle in $G$.
Then one of the following conditions holds:
\item{(i)}the graph $G$ has a planar embedding in which $C$ bounds
a face,
\item{(ii)}there exists a separation $(A,B)$ of $G$ of order at most
three such that $V(C)\subseteq A$ and $G[B]$ does not have a drawing
in a disk with the vertices in $A\cap B$ drawn on the boundary of the disk,
\item{(iii)}there exist two disjoint paths in $G$ with ends
$s_1,t_1\in V(C)$ and $s_2,t_2\in V(C)$, respectively, and otherwise
disjoint from $C$ such that the vertices $s_1,s_2,t_1,t_2$ occur
on $C$ in the order listed.
\smallskip

\proof The lemma is vacuously true for graphs on at most two vertices.
Let $G$ be a graph on at least three vertices, let $C$ be a cycle in $G$,
and assume that the lemma holds for graphs on fewer than $|V(G)|$ vertices.

Suppose first that $G$ is not 3-connected and that  there exists a separation $(A',B')$ of $G$ of order at most two
such that $|A'|,|B'| < |V(G)|$, and assume that the order of  $(A',B')$ is minimum. 
If the order of  $(A',B')$ is two and the two vertices in $A'\cap B'$ are not adjacent, then let $G_1$ be 
obtained from $G[A']$ by adding an edge joining the two vertices in $A'\cap B'$; otherwise let $G_1=G[A']$.
Let $G_2$ be defined analogously, with $A'$ replaced by $B'$.
If $V(C) \subseteq A'$ then $ G_1$ and $C$ satisfy one of  (i),(ii) or (iii) by the choice of $G$. If they satisfy (ii) or (iii)  then the same conclusion is satisfied by $G$ and $C$. 
(If one of the paths as in (iii) uses the added edge, then that edge may be replaced by a path in $G[B']$
joining the  two vertices of $A'\cap B'$. Such path exists by the minimality of $A'\cap B'$.)
If $G_1$ and $C$ satisfy (i) then $G$ and $C$ clearly satisfy (i) or (ii). A symmetric argument applies if 
$V(C) \subseteq B'$. Suppose now that  $V(C) \not \subseteq A'$ and  $V(C) \not \subseteq B'$.
Let $C_1$ and $C_2$ be two cycles obtained from $C[A']$ and $C[B']$, respectively, by adding the edge $e$ joining the two vertices of $A' \cap B'$. 
As before,  if (ii) or (iii) holds for $G_i$ and $C_i$ for some $i \in \{1,2\}$ then the same conclusion holds for $G$ and $C$. Finally, if (i) holds for $G_1$ and $C_1$, and for $G_2$ and $C_2$, then $G$ and $C$ satisfy the same conclusion, as one can combine the embeddings of $G_1$ and $G_2$ by gluing them along~$e$.

Thus we may assume that $G$ is 3-connected.  By [\cite{KelSpecialK4}, Theorem~3.2] either the lemma holds,
or there exists a separation $(A,B)$ of $G$ of order at most three such that
$V(C)\subseteq A$ and $|B-A|\ge 2$.  By moving components of
$G\backslash (A\cap B)$ from $A$ to $B$ we may assume that every component
of $G\backslash  B$ includes at least one vertex of $C$.  We may
assume that $G[B]$ can be drawn in a disk with $A\cap B$ drawn on the
boundary of the disk, for otherwise the lemma holds.  Let $G'$ be obtained
from $G[A]$ by adding an edge joining every pair of nonadjacent vertices in
$A\cap B$.  Then $G'$ satisfies one of (i)--(iii) by the minimality of $G$.
However, since $G[B]$ can be drawn in a disk as specified above,
it follows that $G$ satisfies the same conclusion.\qed

If $G$ is a subdivision of a $3$-connected planar graph, then it has a unique
planar embedding by Whitney's theorem [\cite{WhiCongr}],
and the cycles that bound faces can be
characterized combinatorially.
A cycle $C$ in a graph $G$ is
called \dfn{peripheral} if it is an induced subgraph of $G$,
and $G\backslash V(C)$ is connected.
The following three results are well-known [\cite{TutHowto}, \cite{WhiCongr}].

\newthm{periph}Let $G$ be a subdivision of a $3$-connected planar graph,
and let $C$ be a cycle in $G$. Then the following conditions are
equivalent:
\item{(i)}the cycle $C$ bounds a face in some planar
embedding of $G$,
\item{(ii)}the cycle $C$ bounds a face in every planar
embedding of $G$,
\item{(iii)}the cycle $C$ is peripheral.

\newthm{periphcap}Let $G$ be a subdivision of a $3$-connected planar
graph, and let $C_1,C_2$ be two distinct peripheral cycles in
$G$.  Then the intersection of $C_1$ and $C_2$
 is either null, or a one-vertex graph, or a segment.

\newthm{rotation} Let $G$ be a subdivision of a $3$-connected
planar graph, let $v\in V(G)$
 and let $e_1,e_2,e_3$ be three distinct edges of $G$
incident with $v$.  If there exist peripheral
cycles $C_1,C_2,C_3$ in $G$ such that $e_i\in E(C_j)$
for all distinct indices $i,j\in \{1,2,3\}$, then $v$ has
degree three.

\newsection{disksys}DISK SYSTEMS

The preceding theorems summarize all the properties of peripheral
cycles that we will require.  However, for the sake of greater
generality we will be working with sets of cycles satisfying
only those axioms that will be needed.
Thus we define a \dfn{weak disk system} in a graph $G$ to be
a set $\cal C$ of distinct cycles of $G$, called \dfn{disks}, such that
\item{(X0)}every edge of $G$ belongs to exactly two members of $\cal C$, and
\item{(X1)}the intersection of any two distinct members of $\cal C$
is either null, or a one-vertex graph, or a segment.

The weak disk system under consideration will typically be clear from context, and we will typically refer to elements of a weak disk system $\calc$ simply as disks, rather than disks in $\calc$.

\noindent
A weak disk system is a \dfn{disk system} if it satisfies (X0), (X1) and
\item{(X2)} if $e_1,e_2,e_3$ are three distinct edges incident with a
vertex $v$ of $G$ and there exist disks $C_1,C_2,C_3$ such that
$e_i\in E(C_j)$ for all distinct integers $i,j\in \{1,2,3\}$, then $v$
has degree three.

\noindent
Thus by \refthm{periph}, \refthm{periphcap} and~\refthm{rotation}
the peripheral cycles of a subdivision of a 3-connected planar
graph form a disk system.
If $G'$ is obtained from $G$ by (repeated) rerouting, then
a weak disk system $\calc$ in $G$ induces a weak disk system ${\cal C}'$
in $G'$ in the obvious way. We say that $\calc'$ is the
{\sl weak disk system induced in $G'$ by $\calc$}.
If $\calc$ is a disk system, then so is $\calc'$.

Let $S$ be a subgraph of a graph $H$.
Let us recall that a path $P$ in $H$ is
an \dfn{$S$-path} if it has at least one edge, and
its ends and only its ends belong to $S$.
Now let $\calc$ be a weak disk system in $S$.
An $S$-path $P$ is an $S$-\dfn{jump} if no disk in $\calc$ includes
both ends of $P$.
Let $x_1,x_2,x_3\in V(S)$, let $x\in V(H)-V(S)$,
and let $P_1,P_2,P_3$ be three paths in $H$
such that $P_i$ has ends $x$ and $x_i$, they are pairwise disjoint
except for $x$, and each is disjoint from $V(S)-\{x_1,x_2,x_3\}$.
Assume further that for each pair $x_i,x_j$ there exists a disk
containing both $x_i$ and $x_j$, but no disk contains
all of $x_1, x_2, x_3$.
In those circumstances we say that the triple $P_1,P_2,P_3$ is
an \dfn{$S$-triad}. The vertices $x_1,x_2,x_3$ are its \dfn{feet}. Note that the definition of an $S$-triad depends on   underlying  weak disk system $\calc$. However, we omit  $\calc$ from the notation, as the choice of the weak disk system will be always clear from the context.

Let $S$ be a graph, and let $\calc$ be a weak disk system in $S$.
We say that a subgraph $J$ of $S$ is a \dfn{detached $K_4$-subdivision}
if $J$ is isomorphic to a subdivision of $K_4$,
every segment of $J$ is a segment of $S$, and each of the four
cycles of $J$ consisting of precisely three segments is a disk.


\newthm{claim1}Let $G$ be a graph with no vertices of degree two,
let $S$ be a $G$-subdivision in
a graph $H$,  let $\cal C$ be a weak disk system in $S$, and let $B$ be
an $S$-bridge with at least two attachments such that no disk includes all attachments of $B$.
Then one of the following conditions holds:
\item{(i)}there exists an $S$-jump, or
\item{(ii)}there exists an $S$-triad, or
\item{(iii)}$S$ has a detached $K_4$-subdivision $J$ such that
the attachments of $B$ 
are precisely the branch-vertices of $J$.
\smallskip

\proof
\noindent We may assume that (i) and (ii) do not hold.
Let $S$ and $B$ be as stated, and let $A$ be the
set of all attachments of $B$.
Thus $|A|\ge2$.
Since (i) does not hold, we deduce
that for every pair of elements $a_1,a_2\in A$ there exists a
disk $C\in\cal C$ such that $a_1,a_2\in V(C)$. Since (ii)
does not hold, we deduce that the same holds for every triple of
elements of $A$.

Now let $k\ge3$ be the maximum integer such that for every $k$-element
subset $A'$ of $A$ there exists a disk $C\in\cal C$
such that $A'\subseteq V(C)$. By hypothesis $k<|A|$, and hence
there exist distinct vertices
$a_1,a_2,\ldots,a_{k+1}\in A$
such that $a_1,a_2,\ldots,a_{k+1}\in V(C)$ for no
disk $C\in \cal C$. For $i=1,2,\ldots,k+1$ let
$C_i\in\cal C$ be a disk in $S$ such that $V(C_i)$
includes all of  $a_1,a_2,\ldots,a_{k+1}$ except $a_i$. Then
these disks are pairwise distinct. Since $a_1$ and $a_2$
belong to both $C_3$ and $C_4$ and $\cal C$ satisfies (X1),
there exists
 a segment $P_{12}$ of $S$ that is a subgraph of  $C_3\cap C_4$
and contains $a_1$ and $a_2$.
Similarly,
for all distinct integers $i,j= 1,2,\ldots,k+1$,
there is a segment $P_{ij}$ of $S$ such that
$a_i,a_j\in V(P_{ij})$ and $P_{ij}$ is a subgraph of $C_\ell$ for all
$\ell\in \{1,2,\dots, k+1\}-\{i,j\}$.
Now for all $i=1,2,\dots, k+1$ the vertex
$a_i$ is an end of $P_{ij}$, for otherwise
the segments $P_{ij}$ $(j\in \{1,2,\dots, k+1\}-\{i\})$
would be all equal, implying that
$a_1,a_2,\ldots,a_{k+1}$ all belong to $V(C_t)$ for all
$t=1,2,\ldots,k+1$, a contradiction. Thus $a_1,a_2,\dots, a_{k+1}$
are branch-vertices of $S$.
It follows that $\bigcup P_{ij}$ is a subdivision
of a complete graph $J$. Since $P_{23}\cup P_{24}\cup P_{34}$ is a cycle
and it is a subgraph of $C_1$, it is equal to $C_1$.  Similarly
for $C_2,C_3,C_4$.  Hence $k=3$,
and since (i)  does not hold we deduce from (X1) that $A=\{a_1,a_2,a_3,a_4\}$.
Thus (iii) holds, as desired.~\qed

In the following definitions let
$S$ be a subgraph of a graph $H$ and let
$\cal C$ be a weak disk system in
a graph $S$.
Let $C\in\cal C$, and let $P_1$ and $P_2$ be
two disjoint $S$-paths with ends $u_1,v_1$ and $u_2,v_2$,
respectively, such that $u_1,u_2,v_1,v_2$ belong to $V(C)$ and
occur on $C$ in the order listed. In those circumstances we say
that the pair $P_1,P_2$ is an \dfn{$S$-cross}. We also say that
it is an \dfn{$S$-cross on $C$.} We say that $u_1,v_1,u_2,v_2$
are the \dfn{feet} of the cross. We say that the cross
$P_1,P_2$ is \dfn{weakly free} if
\item{(F1)}for $i=1,2$ no segment of $S$ includes both ends of $P_i$.

\noindent
We say that a cross $P_1,P_2$ is \dfn{free} if it satisfies (F1) and
\item{(F2)}
no two segments of $S$ that share a vertex include all the feet of the cross.


The intent of freedom is that the feet of the cross are not separated from ``most of $S$"
by a separation of order at most three, but it does not quite work that way for
our definition. If $C$ is a cycle in $S$ consisting of three segments, then
no free cross on $C$ has that property.
That should be regarded as a drawback of our definition.
However, it turns out that it is not a problem in any of our applications,
because in all applications the graph $G$ has girth at least four.
On the other hand, there does not seem to be an easy way to eliminate
crosses on cycles consisting of three segments, and since we do not
need to do it, we chose to avoid it.
It should be noted, however, that the ``right" definition of freedom should avoid
crosses on cycles consisting of three segments.

\medskip
A separation $(X,Y)$ of $H$ is called
an $S$-\dfn{separation} if the order of $(X,Y)$ is at most three, $X-Y$
includes at most one branch-vertex of $S$,
and the
graph $H[X]$ does not have a drawing in a disk with $X\cap Y$ drawn on the
boundary of the disk.

We say that ${\cal C}$
is \dfn{locally planar} in $H$ if for every $S$-bridge $B$ of $H$
with at least two attachments there exists
a disk $C_B \in \calc$
such that $C_B$ includes all attachments of $B$
and for every disk $C\in\calc$ the graph $C\cup \bigcup B$ has a planar drawing
with $C$ bounding the outer face, where the union is taken over
all $S$-bridges $B$ of $H$ with $C_B=C$.

Let $Z$ be a segment of $S$, let $z,w$ be the ends of $Z$,
and let $P_1,P_2$ be two disjoint $S$-paths in $H$
with ends $x_1,y_1$ and $x_2,y_2$, respectively, such that
$z,x_1,x_2,y_1,w\in V(Z)$ occur on $Z$ in the order listed,
and $y_2\not\in V(Z)$. Let $P_3$ be a path disjoint from
$V(S)-\{y_2\}$ with one end  $x_3\in V(P_1)$
and the other  $y_3\in V(P_2)$ and otherwise disjoint from
$P_1\cup P_2$.
We say that the triple $P_1,P_2,P_3$ is an \dfn{$S$-tripod}
\dfn{based at} $Z$, and that $x_1,y_1,x_2,y_2$ (in that order) are its
\dfn{feet}.
We say that $zZx_1$, $y_1Zw$ and $y_3P_2y_2$ are the {\sl legs} of the tripod.
See Figure~\reffig{tripod}.

\bigskip

\bigskip
\centerline{\includegraphics[scale=1.70]{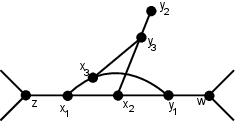}}
\bigskip
\centerline{Figure~\newfig{tripod}. An $S$-tripod.}
\bigskip




\newthm{planarlemma}Let $G$ be a  graph with
no vertices of degree two,
let $H$ be
a connected graph, and let $S$ be a $G$-subdivision
in $H$ with a weak disk system ${\cal C}$.
Then $H$ has a $G$-subdivision $S'$ obtained from $S$ by repeated
I-reroutings such that $S'$ and the weak disk system ${\cal C}'$
in $S'$ induced by ${\cal C}$ satisfy
one of the following conditions:
\item{(i)}there exists an $S'$-jump, or
\item{(ii)}there exists a weakly free $S'$-cross in $H$
on some member of ${\cal C}'$, or
\item{(iii)}$H$ has an $S'$-separation $(X,Y)$ such that $X-Y$ includes
no branch-vertex of $S'$, or
\item{(iv)}$S'$ has a detached $K_4$-subdivision $J$ and
$H$ has an $S'$-bridge $B$ such that the attachments
of $B$ are precisely the branch-vertices of $J$, or
\item{(v)}there exists an $S'$-triad, or
\item{(vi)} the weak disk system ${\cal C}'$ is locally
planar in $H$.
\smallskip

\proof
We proceed by induction on $|V(H)|$.
Suppose for a contradiction that none of (i)--(vi) holds.
We start with the following claim.

\claim{1} Let $S'$ be a $G$-subdivision in $H$ obtained from $S$ by
  repeated I-reroutings, and let ${\cal C}'$ be the weak disk system
  induced in $S'$ by $\cal C$. Then for every $S'$-bridge $B$ of $H$
  with at least two attachments
  there exists a disk $C\in{\cal C}'$ such that $V(C)$ includes all
  attachments of $B$.

\noindent
Claim (1) follows from~\refthm{claim1}, for otherwise one of the outcomes
(i), (iv), (v) holds, a contradiction.
This proves (1).

\claim{2} There exists a $G$-subdivision $S'$ in $H$ obtained from $S$ by
  repeated I-reroutings such that every $S'$-bridge is \stable.

\noindent
To prove (2) let $S'$ be as in~\refthm{stable}.
We may assume that there exists an unstable $S'$-bridge $B'$, for otherwise
(2) holds.
Let $Z$ be a segment of $S$ that includes all attachments of $B'$.
By~\refthm{stable} there exist a separation $(X,Y)$ and vertices
$x,y\in V(Z)$ such that $Y$ includes every branch-vertex of $S$,
$V(B'\cup xZy)\subseteq X$
and $X\cap Y=\{x,y\}$.
Since $(X,Y)$ does not satisfy (iii), the graph $H[X]$
has a drawing in a disk with $x,y$ on the boundary of the disk.
Let $H'$ be obtained from $H\backslash (X-Y)$ by adding an edge
joining $x,y$ if $x$ and $y$ are distinct and not adjacent in $H$,
and let $H':=H\backslash (X-Y)$ otherwise.
By induction applied to $G$, $H'$ and a suitable modification of the graph $S$
we conclude that $H'$ satisfies one of the conclusions of the lemma.
But then $H$ also satisfies the conclusion of the lemma, because
if $H'$ satisfies (vi) it follows from the planarity of $H[X]$ that
so does $H$.
The other conditions are straightforward.
This proves (2).

\claim{3} There exist a $G$-subdivision $S'$ in $H$ obtained from $S$ by
  repeated I-reroutings and an $S'$-tripod.

\noindent
To prove (3) we choose $S'$ as in (2); hence every
$S'$-bridge is \stable.
It now follows from (1) that
for every $S'$-bridge $B$ there exists a unique disk $C$ in the weak
disk system ${\cal C}'$ induced in $S'$ by $\cal C$
such that all attachments of $B$ belong to $V(C)$.
 For every disk $C$ of $S'$ let $H_C$ be the union
of $C$ and all $S'$-bridges $B$ whose attachments are included
in $V(C)$. Since (vi) does not hold,
there exists a disk $C$ of $G$ such that $H_C$
does not have a planar drawing with $C$ bounding the infinite
face.

By \refthm{2paths}
and the fact that (iii) does not hold there
exists an $S'$-cross $P_1,P_2$ in $C$. For $i=1,2$, let $x_i,y_i$ be the ends
of $P_i$ and let $B_i$ be the $S'$-bridge that includes $P_i$.
We may assume that there is a segment $Z$ of $S'$ such that
$x_1,x_2,y_1\in V(Z)$, for otherwise the $S'$-cross $P_1,P_2$ satisfies (F1).
We claim that we may assume that $y_2\not\in V(Z)$.
Indeed, if $y_2\in V(Z)$, then since $B_1$ and $B_2$ are \stable,
there exists a path from $P_1\cup P_2$, say from $P_2$, to
a vertex $v\in
V(C)-V(Z)$, disjoint from $V(P_1\cup P_2\cup S')-\{v\}$. It follows
that $P_1\cup P_2\cup P$ includes a $S'$-cross with at least one
foot outside $Z$.
Thus we may assume that $y_2\notin V(Z)$.

If $B_1=B_2$, then there exists a path $P_3$ as in the definition of
$S'$-tripod, and hence the claim holds. Thus we may assume that $B_1\not=B_2$.
Since $B_1$ is \stable\ there exists a path $P_3$ in $B_1$ with one
end in $V(P_1)-\{x_1,y_1\}$ and the other end  $z\in V(C)-V(Z)$.
If $z=y_2$, then $P_1,P_2,P_3$ is an $S$-tripod,
 as desired, and so we may assume that $z\ne y_2$.
Then $P_1\cup P_2\cup P_3$ includes a weakly free $S'$-cross, unless
some segment $Z'$ of $S'$ includes either $z,y_2,y_1$ in the order listed,
or $z,y_2,x_1$ in the order listed.
By symmetry we may assume the former.
Then $y_1$ is a common end of $Z$ and $Z'$.
Let $S''$ be obtained from $S'$ by replacing $x_1Zy_1$ by $P_1$;
then $P_3,P_2\cup x_1Zx_2$ is a weakly free $S''$-cross, as desired.
This proves (3).
\medbreak

To complete the proof of the theorem
we may select a $G$-subdivision $R$ in $H$ obtained from $S$ by
repeated I-reroutings and an $R$-tripod $P_1,P_2,P_3$
such that the sum of the lengths of the tripod's legs is minimum.
Let $Z,z,w,x_1,y_1,x_2,y_2,x_3,y_3$ be as in the definition of tripod.

Let $R'$ be obtained from $R$ by rerouting $x_1Zy_1$
along $P_1$; then $x_1Zy_1$, $P_3\cup y_3P_2y_2$,
$x_2P_2y_3$ is an $R'$-tripod with the
same legs. Thus there is symmetry between
$x_1Zy_1\cup x_2P_2y_3$
and $P_1\cup P_3$.

Let $X'$ be the vertex-set of
$x_1Zy_1\cup x_2P_2y_3
\cup P_1\cup P_3$, and let
$Y'=V(R)-(X'-\{x_1,y_1,y_3\})$. If there is no path between
$X'$ and $Y'$ in $H\backslash \{x_1,y_1,y_3\}$, then there exists
a separation $(X,Y)$ of order three in $H$ with $X'\subseteq X$
and $Y'\subseteq Y$ (and hence $X\cap Y=\{x_1,y_1,y_3\}$). Then
$(X,Y)$ is an $R$-separation, and hence (iii) holds, a contradiction.
Thus there exists a path
$P$ in $H\backslash \{x_1,y_1,y_3\}$ with ends $x\in X'$ and
$y\in Y'$. From the symmetry established in the previous paragraph
we may assume that $x\in V(P_1)\cup V(P_3)-\{x_1,y_1,y_3\}$.
It follows from the minimality of legs that $y\not\in V(Z)\cup V(P_2)$.

Let $C_1,C_2$ be the two disks of $R$ that include  $Z$.
Then $y_2\in V(C_i)$ for some $i=1,2$, say $i=1$,
for otherwise $P_2$ is an $R$-path satisfying (i).
Since $y_2\not\in V (Z)$, (X1) implies that
$y_2\not\in V(C_2)$.
Since the vertices
$x_1,y_1,y_2,y$ are attachments of an $R$-bridge, by (1) there exists a
disk $C$ in $G$ such that $x_1,y_1,y_2,y\in V(C)$.
Since $x_1,y_1\in V(C)$, (X1) implies that $C=C_1$ or $C=C_2$, but
$y_2\not\in V(C_2)$, and so $C=C_1$.  In particular, $y,y_2\in V(C_1)$.
Since $y\ne y_2$ (because $y\not\in V(P_2)$),
 $P_1\cup P_2\cup P_3\cup P$ includes an $R$-cross in $C_1$
satisfying (F1),
unless either $z=x_1$ and $z,y_2,y$ appear on a segment incident with $z$ 
in the order listed,
or $y_1=w$ and $w,y_2,y$ appear on a segment incident with $w$ in the order listed.
We may therefore assume by symmetry that the former case holds.
Let $R''$ be obtained from $R$ by replacing $x_1Zy_1$ by $P_1$;
then $y_1Zx_2\cup P_2,P\cup P_3$ includes  an $R''$-cross satisfying (F1),
as desired.\qed

\medbreak
Our next objective is to improve outcome (ii) of the previous lemma.
Let
$S$ be a subgraph of $H$, let $C$ be a cycle in $S$, and let
$P_1,P_2$ be a weakly free $S$-cross on $C$.  If the cross
$P_1,P_2$ is not free, then there exist two distinct segments
$Z_1,Z_2$ of $S$, both subgraphs of $C$ and both incident with
a branch-vertex $v$ of $S$ such that
$Z_1\cup Z_2$ includes all the feet of $P_1,P_2$.
In those circumstances we say that the cross $P_1,P_2$ is \dfn{centered at}
$v$ and that it is \dfn{based at} $Z_1$ and $Z_2$.
We will treat the cases when $v$ has degree three and when it has
degree at least four separately.

We say that an $S$-triad in a graph $H$ is \dfn{local} if there exists a
vertex $v$ of $S$ of degree three in $S$ such that each of the three
segments of $S$ incident with $v$ includes exactly one foot of the
triad.
We say that the local
$S$-triad is \dfn{centered} at $v$.

\newthm{weakfree} Let $G$ be a graph with no vertices of degree two,
let $H$ be a graph, let $S$ be a $G$-subdivision in $H$ with a
weak disk system $\calc$,
let $C\in\calc$, let $v\in V(C)$ have degree in $S$ exactly three,
 and let $P_1,P_2$ be a weakly free
$S$-cross in $H$ on $C$ centered at $v$.
Then there exist a $G$-subdivision $S'$ obtained
from $S'$ by exactly one {\rm T}-rerouting
centered at $v$ and a local $S'$-triad.

\proof
For $i=1,2$ let $x_i,y_i$ be the ends of $P_i$, and let
$P_1,P_2$  be based at $Z_1$ and $Z_2$.
Then we may assume that $x_1,x_2,v\in V(Z_1)$ occur on
$Z_1$ in the order listed; then $y_2,y_1,v\in V(Z_2)$
occur on $Z_2$ in the order listed.  Let $S'$ be the
$G$-subdivision obtained from $S$ by rerouting $vZ_2y_2$ along $P_2$.
Then $P_1,y_1Z_2y_2, vZ_2y_1$ is a desired $S'$-triad.~\qed

Converting weakly free crosses centered at vertices of degree
at least four into free crosses is best done by
splitting vertices, but we are concerned with subdivisions, and
therefore we take a different route.
In the next lemma we need $\cal C$ to be a disk system
(not merely a weak one).

\newthm{weakfree2}Let $G$ be a graph with no vertices of degree
two, let $H$ be a graph, let $S$ be a $G$-subdivision in $H$ with a
disk system $\calc$,
and assume that $H$ has a weakly free $S$-cross
centered at a vertex of degree at least four.
Then there exists a $G$-subdivision $S'$ obtained from $S$ by
repeated rerouting  such that $S'$ and the disk system ${\cal C}'$
induced in $S'$ by $\cal C$ satisfy one the following conditions:
\item{(i)} $H$ has an $S'$-jump,
\item{(ii)} $H$ has a free $S'$-cross on some disk in ${\cal C}'$, or
\item{(iii)} $H$ has an $S'$-separation $(X,Y)$ such that $X-Y$ includes
no branch-vertex of $S'$.
\smallskip

\proof
Let $P_1,P_2$ be a weakly free $S$-cross in $H$
centered at a vertex $v$ of degree at least four.
Thus there exist
two segments $Z_1,Z_2$ of $S$, both incident with $v$, such that
 $Z_1,Z_2$ include all the feet of the cross.
For $i=1,2$ let $x_i,y_i$ be the ends of $P_i$.
We may assume that $x_1,x_2,v\in V(Z_1)$  occur on $Z_1$ in the order
listed; then $y_2,y_1,v\in V(Z_2)$ and they occur on $Z_2$ in the
order listed.
For $i=1,2$ let $v_i$ be the other end of $Z_i$ and let
$L_1=x_1Z_1v_1$ and $L_2=y_2Z_2v_2$.

Consider all triples $(S', P_1',P_2')$, where $S'$ is
a $G$-subdivision  obtained from $S$ by repeated rerouting
and $P_1',P_2'$ is a weakly free $S'$-cross based at $Z_1',Z_2'$
(where $Z'_1,Z'_2$ are the branches of $S'$ corresponding to $Z_1,Z_2$).
We may assume that among all such triples 
the triple $(S,P_1,P_2)$ is chosen so that
$|V(L_1)|+|V(L_2)|$ is minimum.

Let $X'$ be the vertex-set of
$P_1\cup P_2\cup vZ_1x_1\cup vZ_2y_2$
and let $Y'=V(S)-(X'-\{v,x_1,y_2\})$. If there is no path
in $H\backslash\{v,x_1,y_2\}$ with one end in $X'$ and the other
in $Y'$, then there exists a separation $(X,Y)$ of order three with
$X'\subseteq X$ and $Y'\subseteq Y$. This separation satisfies (iii),
and so we may assume that  there exists a path $P$ in
$H\backslash\{v,x_1,y_2\}$
with one end $x\in X'$ and the other end $y\in Y'$. From the symmetry
we may assume that $x$ belongs to the vertex-set of
$P_1\cup vZ_2y_2$.

If $y\in V(L_1)$, then replacing $P_1$ by $P$ if $x\not\in V(P_1)$
and by $P\cup xP_1y_1$ otherwise produces a cross that contradicts
the choice of the triple $(S,P_1,P_2)$.
If $y\in V(L_2)$, then replacing $yZ_2x$ by $P$ if $x\not\in V(P_1)$
and replacing $yZ_2y_1$ by $P\cup xP_1y_1$ results in a $G$-subdivision
$S'$ obtained from $S$ by a rerouting, and $P_1,P_2$ can be modified
to give a cross $P_1',P_2'$ such that the triple $(S',P_1',P_2')$
contradicts the choice of $(S,P_1,P_2)$. Thus $y\not\in V(Z_1\cup Z_2)$.

Let $C$ be the disk that includes both $Z_1$ and $Z_2$
(it exists, because $P_1,P_2$ is a cross), and for $i=1,2$ let
$C_i$ be the other disk that includes $Z_i$.
If $y\in V(C)$, then $P_1\cup P_2\cup P$ includes a free cross,
and so (ii) holds. Thus we may assume that $y\not\in V(C)$.
Similarly, if $y\not\in V(C_2)$, then $P_1\cup P$ includes an $S$-jump
with one end $y$ and the other end $x$ or $y_1$, and so we may assume that
$y\in V(C_2)$.
Since  $v$ has degree at least four,
(X1) and (X2) imply that $V(C_1)\cap V(C_2)=\{v\}$.
It follows that $y\not\in V(C_1)$.
Now let $S'$ be obtained from $S$ by an X-rerouting using the cross
$P_1,P_2$, and let $Z_1',Z_2'$ be the segments of $S'$ corresponding
to $Z_1,Z_2$, respectively.
Thus $Z_1'=v_1Z_1x_1\cup P_1\cup y_1Z_2v$ and
$Z_2'=v_2Z_2y_2\cup P_2\cup x_2Z_1v$.
Now $P\cup xZ_2y_1$ includes  an $S'$-jump with one end $y$ and the
other end in the interior of $Z_1'$, and so (i) holds.~\qed

We can summarize some of the lemmas of this section as follows.

\newthm{summary1}Let $G$ be a connected graph with no vertices of degree two
that is not the complete graph on four vertices,
let $H$ be a graph, and let $S$ be a $G$-subdivision
in $H$ with a disk system ${\cal C}$. Then
$H$ has a $G$-subdivision $S'$ obtained from $S$ by repeated
reroutings such that $S'$ and the weak disk system ${\cal C}'$
induced in $S'$ by $\cal C$ satisfy
one of the following conditions:
\item{(i)}there exists an $S'$-jump in $H$, or
\item{(ii)}there exists a free $S'$-cross in $H$
on some disk of ${\cal C}'$, or
\item{(iii)}$H$ has an $S'$-separation $(X,Y)$ such that $X-Y$ includes
no branch-vertex of $S'$, or
\item{(iv)}there exists an $S'$-triad, or
\item{(v)} the disk system ${\cal C}'$ is locally
planar in $H$.

\proof
By~\refthm{planarlemma} there exists a $G$-subdivision $S_1$ obtained from
$S$ by a sequence of reroutings such that one of the outcomes of
that lemma holds.
But~\refthm{planarlemma}(iv) does not hold, because $\cal C$ satisfies (X2)
and $G$ is not $K_4$.
We may assume therefore that \refthm{planarlemma}(ii) holds, for otherwise
$S_1$ and the weak disk system induced in $S_1$ by $\cal C$ satisfy
\refthm{summary1}.
Thus $S_1$ has a disk $C$ and a weakly free cross $P_1,P_2$ on $C$.
We may assume that $P_1,P_2$ is not free, for otherwise \refthm{summary1}(ii)
holds.
Thus there exists a branch-vertex $v$ of $S_1$ that belongs to $C$
and two distinct segments
$Z_1,Z_2$ of $S_1$, both subgraphs of $C$ and both incident with $v$
such that the cross $P_1,P_2$ is centered at $v$ and based at $Z_1,Z_2$.
If $v$ has degree three in $S_1$, then the lemma holds by~\refthm{weakfree}
and if $v$ has degree at least four,
then the lemma holds by~\refthm{weakfree2}.~\qed

The following theorem will be used in~[\cite{HegTho}].
Recall that a graph $G$ is
\afc\ if it is 3-connected, has at least five vertices,
 and, for every separation
$(A,B)$ of $G$ of order 3, one of $A-B, B-A$ contains at most one vertex.

\newthm{rajneesh}Let $G$ be a  graph 
with no vertices of degree two that is not the complete graph on four vertices,
let $H$ be an \afc\ graph, and let $S$ be a $G$-subdivision
in $H$ with a disk system ${\cal C}$. Then
$H$ has a $G$-subdivision $S'$ obtained from $S$ by repeated
reroutings such that $S'$ and the disk system ${\cal C}'$
induced in $S'$ by $\cal C$ satisfy
one of the following conditions:
\item{(i)}there exists an $S'$-jump in $H$, or
\item{(ii)}there exists a free $S'$-cross in $H$
on some disk of $S'$, or
\item{(iii)}there exists an $S'$-triad, or
\item{(iv)} the disk system ${\cal C}'$ is locally
planar in $H$.

\proof
By~\refthm{summary1} we may assume that there exists a $G$-subdivision
$S'$ obtained from $S$ by repeated
reroutings such that $S'$ and the weak disk system ${\cal C}'$
induced in $S'$ by $\cal C$ satisfy~\refthm{summary1}(iii),
for otherwise the theorem holds.
Thus $H$ has an $S'$-separation $(X,Y)$ such that $X-Y$ includes
no branch-vertex of $S'$.
Since $S$ has at least five branch-vertices, it follows that
$|Y-X|\ge2$. But $|X-Y|\ge2$, because $H[X]$ cannot be drawn in
a disk with $X\cap Y$ drawn on the boundary of the disk.
This contradicts the \afcity\ of $H$.~\qed

\newsection{triads}TRIADS

In this section we improve outcome (iv) of \refthm{summary1}.
A graph $G$ is {\sl\ifc} if it is 3-connected and
for every separation $(A,B)$ of order three one of $G[A], G[B]$ has
at most three edges.
(Thus  every 4-connected graph is \ifc\ and every \ifc\ graph is
\afc.)
If $S$ is a $G$-subdivision in a graph $H$, then there is a mapping
$\eta$ that assigns to each $v\in V(G)$ the corresponding vertex
$\eta(v)\in V(S)$, and to every edge $e\in E(G)$ the corresponding
path $\eta(e)$ of $S$.
We say that $\eta$ is a \dfn{\he}, and we write $\eta:G\emb S\subseteq H$
to denote the fact that $\eta$ is a \he\ that maps $G$ onto the
subgraph $S$ of $H$.


\newthm{lemma2}Let $G$ be an \afc\ graph,
let $H$ be a graph, let $S$ be a $G$-subdivision in $H$ with a
weak disk system ${\cal C}$,
and assume that there exists an $S$-triad in $H$ that is not local
and has set of feet $F$.
Assume also that if $|V(G)|\le 6$, then $G$ is \ifc.
Then
\item{(i)}there exists a segment of $S$ with both ends in $F$, or
\item{(ii)}$S\backslash F$ is connected and $F$ is an
independent set in $S$.
\smallskip

\proof
Let the $S$-triad be
$Q_1,Q_2,Q_3$, and let
$F=\{x_1,x_2,x_3\}$ be  labeled so that $x_i$ is an end of $Q_i$.
Let $J$ be the subgraph of $S$ with vertex-set $F$ and no edges.
We may assume that $S$ has at least two $J$-bridges, for otherwise
(ii) holds.  Assume first that some $J$-bridge of $S$ includes no
branch-vertex of $S$, except possibly as an attachment.
Then that $J$-bridge is a subgraph of
 a segment $Z$ that includes two members of $F$, say $x_1$
and $x_2$.  It follows that $x_1$ and $x_2$ are the ends of $Z$, for
if $x_1$ is an internal vertex of $Z$, then the disk containing
$x_1$ and $x_3$ contains $x_2$ as well, a contradiction.  Hence (i) holds.

Thus we may assume that
every $J$-bridge of $S$ includes a branch-vertex of $S$ that is not
an attachment of $J$, and since
there are at least two $J$-bridges of $S$, it follows that
$S$ has a separation $(X,Y)$ with $X\cap Y=\{x_1,x_2,x_3\}$ such that
both $X-Y$ and $Y-X$ include a branch-vertex of $S$.

We claim that one of $X-Y, Y-X$ includes at most one branch-vertex of $S$.
To prove this claim, suppose the contrary and let
$\eta:G\hookrightarrow S\subseteq H$ be a \he.
Let $z_1,z_2,z_3 \in V(G)\cup E(G)$ be defined as follows.
Let $i\in \{1,2,3\}$.  If $x_i$ is a branch-vertex of $S$, then let
$z_i\in V(G)$ be such that $\eta (z_i)=x_i$; otherwise $x_i$ is the
interior vertex of a unique segment $\eta (z_i)$
of $S$, and we define $z_i$ that way.
Let $X'$ be the set of all vertices $x$ of $G$ such that
$\eta (x)\in X$,
and let $Y'$ be defined
analogously.  Then $X'\cup Y'=V(G)$, and there are exactly $3-|X'\cap Y'|$
edges of $G$ with one end in $X'-Y'$ and the other in $Y'-X'$. Note
that $X'\cap Y'=\{z_1,z_2,z_3\}\cap V(G)$.  If $z_1,z_2,z_3\in V(G)$,
then our claim follows from the \afcity\ of $G$.  For the next case
assume that $z_1\in E(G)$ and $z_2,z_3\in V(G)$, and let $u_1,v_1$ be the
ends of $z_1$ with $u_1\in X'$ and $v_1\in Y'$.  By the \afcity\ of $G$
applied to the separation $(X'\cup\{v_1\},Y')$ we deduce that
$|Y'-X'-\{v_1\}|\le 1$, and, by symmetry, $|X'-Y'-\{u_1\}|\le 1$.
Thus $|V(G)|\le 6$, and hence $G$ is \ifc.
Since $G$ has at least five vertices we may assume that
$X'-Y'-\{u_1\}$ is not empty, say $x\in X'-Y'-\{u_1\}$.
Then $x$ has neighbors $u_1,z_2,z_3$.  Since $u_1$ has degree at least three,
it is adjacent to $z_2$ or $z_3$, and hence $x$ has degree three and
belongs to a triangle, contrary to the \ifcity\ of $G$.
Thus, the claim holds when
at most one of $z_1,z_2,z_3$ is an edge.  The other two cases are
similar, and are omitted.  This completes the proof of our claim that
one of $X-Y, Y-X$ includes exactly one branch-vertex of $G$.
  From the symmetry we may assume that $X-Y$ includes exactly one
branch-vertex of $S$, say $v$.  It follows that $v$ has degree three
and that $Q_1,Q_2,Q_3$ is a local triad, a contradiction.~\qed

If $G$ is \ifc\ and planar we have the following corollary.

\newthm{planartriad}Let $G$ be an \ifc\ planar graph,
let $H$ be a graph,  let $S$ be a $G$-subdivision in $H$,
and let $\cal C$ be the disk system in $S$ consisting of peripheral
cycles of $S$.
Then every $S$-triad is local.

\proof
Let $F$ be the set of feet of an $S$-triad, and let us assume for a contradiction
that the $S$-triad is not local.
Let us fix a drawing of $S$ in the sphere.
Since every pair of vertices in $F$ are cofacial by~\refthm{periph},
there exists a simple closed curve $\phi$ intersecting $S$ precisely in
the set $F$ and such that both disks bounded by $\phi$ include
a branch-vertex of $S$.
Thus~\refthm{lemma2}(ii) does not hold, and~\refthm{lemma2}(i) does
not hold by the \ifcity\ of $G$ and
the fact that the $S$-triad is not local.
That contradicts~\refthm{lemma2}.~\qed






\newthm{localtriadiv} Let $G,H$ be graphs, where $G$ has no vertices
of degree two, let $S$ be a $G$-subdivision
in $H$, let $\calc$ be a weak disk system in $S$, and
let $Q_1,Q_2,Q_3$ be an $S$-triad in $H$
such that two of its feet are the ends of a segment $Z$ of $S$. Then
there exist a $G$-subdivision $S'$ in $H$
obtained from $S$ by {\rm I}-rerouting the segment $Z$ and an $S'$-jump.

\proof Let the feet of the $S$-triad be $x_1,x_2,x_3$ be numbered so that
$x_i$ is an end of
$Q_i$, and let $Z$ have ends $x_1$ and $x_2$.  Let $S'$ be obtained
from $S$ by replacing $Z$ by $Q_1\cup Q_2$.  Then $Q_3$ is an
$S'$-jump, for if its ends belong to a disk of $S'$, then the
corresponding disk of $S$ would include all of $x_1,x_2,x_3$, contrary to the
definition of an $S$-triad. This proves the lemma.~\qed

Let $G,H$ be graphs, let $G$ have no vertices of degree two, let
$S$ be a $G$-subdivision in $H$,
let $v\in V(S)$ have degree three in $S$, let $Z_1,Z_2,Z_3$ be the
three segments of $S$ incident with $v$, and let $Q_1,Q_2,Q_3$ be a local
$S$-triad centered at $v$ with feet $x_1,x_2,x_3$, where $x_i\in V(Z_i)$.
Let $S'$ be the $G$-subdivision obtained from $S\cup Q_1\cup Q_2
\cup Q_3$ by deleting $v$ and all the edges and internal vertices
of the paths $x_iZ_iv$ for $i=1,2,3$.
We say that $S'$ was obtained from $S$ by a \dfn{triad exchange}.
It follows that $x_1Z_1v, x_2Z_2v, x_3Z_3v$ is an $S'$-triad.

\newthm{localtriad} Let $G,H$ be graphs, where $G$ has no vertices
of degree two, let $S$ be a $G$-subdivision
in $H$,
let $\calc$ be a weak disk system in $S$, and let $Q_1,Q_2,Q_3$ be a local
$S$-triad in $H$.
Then there exists a $G$-subdivision $S'$
obtained from $S$ by repeated rerouting
and at most one triad exchange
such that $S'$ and the weak disk system ${\cal C}'$ in $S'$
induced by $\cal C$ satisfy one of the following conditions:
\item{(i)} there exists an $S'$-jump in $H$, or
\item{(ii)} there exists a free $S'$-cross on some member of ${\cal C}'$, or
\item{(iii)} there exists an $S'$-separation in $H$.
\smallskip

\proof
Let the triad $Q_1,Q_2,Q_3$ be centered at $v$,
let its feet be $x_1,x_2,x_3$, let
$Z_1,Z_2,Z_3$ be the three segments of $S$ incident with $v$ numbered
so that  $x_i\in V(Z_i)$,
 and let $v_i$ be the other end of $Z_i$.
Let $L_i$ be the subpath of $Z_i$ with ends $v_i$ and $x_i$,
and let $P_i$ be the subpath of $Z_i$ with ends $v$ and $x_i$.
We say that the paths $L_1$, $L_2$, $L_3$ are the \dfn{legs} of the $S$-triad
$Q_1$, $Q_2$, $Q_3$.
We may assume that $S$ and $Q_1,Q_2,Q_3$ are chosen so that
there is no $G$-subdivision of $H$
with a triad as above
obtained from $S$ by a  rerouting
such that the sum of the lengths of its
legs is strictly smaller than $|E(L_1)|+|E(L_2)| + |E(L_3)|$.
Let $X_1=V(P_1\cup P_2\cup P_3\cup Q_1\cup Q_2\cup Q_3)$ and
$Y_1=V(S)-(X_1-\{x_1,x_2,x_3\})$.
If $H\backslash\{x_1,x_2,x_3\}$ has
no path between $X_1$
and $Y_1$, then
$H$ has a separation $(X,Y)$ such that $X\cap Y=\{x_1,x_2,x_3\}$,
$X_1\subseteq X$, and $Y_1 \subseteq Y$.
Then $(X,Y)$ satisfies (iii), as desired.

We may therefore assume that there exists a path $P$ as above. Let
the ends of $P$ be $x\in X_1- \{x_1,x_2,x_3\}$
and $y\in Y_1-\{x_1,x_2,x_3\}$.
We may assume that $P$ has no internal vertex in $X_1\cup Y_1$.
We claim that $y\notin V(L_1\cup L_2\cup L_3)$.
Indeed, if $x\in V(Q_1\cup Q_2\cup Q_3)$, then $y\notin V(L_1\cup L_2\cup
L_3)$ by the choice of $Q_1,Q_2,Q_3$ (no change of $S$ needed).
So we may assume that $x\in V(P_1\cup P_2\cup P_3)$ and
 $y\in V(L_1\cup L_2\cup L_3)$.
Now replacing a path
of $P_1\cup P_2\cup P_3$ by $P$ is an I-rerouting or  T-rerouting, and the
resulting $G$-subdivision $S'$ has an $S'$-triad that contradicts the choice of
$S$ and $Q_1,Q_2,Q_3$.  Thus $y\notin V(Z_1\cup Z_2\cup Z_3)$.

The operation of triad exchange exchanges the roles of
$P_1\cup P_2\cup P_3$ and $Q_1\cup Q_2\cup Q_3$.
Thus by applying the triad exchange operation if needed we gain
symmetry between $P_1\cup P_2\cup P_3$ and $Q_1\cup Q_2\cup Q_3$.
Thus may assume that $x\in V(P_1\cup P_2\cup P_3)$.
We may assume that $S$ and $P$ do not satisfy
(i), and hence  there exists a disk
$C$ in $S$ such that $x,y\in V(C)$. It follows that
$C$ includes two of the segments incident with $v$, say $Z_1$ and $Z_2$.
We may assume that
$Q_1\cup Q_2, P$ is not a free $S$-cross
in $C$ for otherwise (ii) holds, and
hence $v_1=x_1$, $v_2=x_2$ and there is a
segment $Z$  of $S$ with ends $v_1$ and $v_2$.
Let $S'$ be the $G$-subdivision obtained from $S$ by the triad exchange
that replaces $Q_1,Q_2,Q_3$ by $P_1,P_2,P_3$.
Then $P\cup P_1\cup P_2\cup P_3$ includes an $S'$-path with ends
$x_3$ and $y$. We may assume that this path is not an $S'$-jump,
for otherwise (i) holds.
Thus there exists a disk $C'$ in $S'$ that includes $Z$ and $x_3$,
and hence includes all of $x_1,x_2,x_3$, contrary to the fact that
$Q_1,Q_2,Q_3$ is a triad.~\qed


The results
thus far can be summarized as follows.

\newthm{summary2}Let $G$ be an \afc\  graph, let $H$ be
a graph, and let $S$ be a $G$-subdivision
in $H$ with a disk system ${\cal C}$.
Assume that if $G$ has at most six vertices, then it is \ifc.
Then $H$ has a $G$-subdivision $S'$ obtained from $S$ by repeated
reroutings and possibly one triad exchange
such that $S'$ and the disk system ${\cal C}'$
induced in $S'$ by $\cal C$ satisfy
one of the following conditions:
\item{(i)}there exists an $S'$-jump in $H$, or
\item{(ii)}there exists a free $S'$-cross in $H$
on some disk of $S'$, or
\item{(iii)}$H$ has an $S'$-separation, or
\item{(iv)}there exists an $S'$-triad with set of feet $F$ such that
$S'\backslash F$ is connected and $F$ is an
independent set in $S'$, or
\item{(v)} the disk system ${\cal C}'$ is locally
planar in $H$.

\proof
Let $S'$ be as in~\refthm{summary1}, and let $\calc'$ be the
corresponding disk system in $S'$.
We may assume that~\refthm{summary1}(iv) holds, for otherwise
the lemma holds.
Let $t$ be an $S'$-triad.
If $t$ is local, then the result holds by~\refthm{localtriad}.
Otherwise by~\refthm{lemma2}
and~\refthm{localtriadiv} either outcome (i) or outcome (iv) holds.~\qed

\newsection{planar}WHEN $G$ IS PLANAR




We are now ready to reformulate the above results in terms of embedded
graphs.  By a \dfn{surface} we mean a compact connected 2-dimensional
manifold with no boundary.  A graph $S$ embedded in a surface
$\Sigma$ is \dfn{polyhedrally embedded} if $S$ is a subdivision of a
$3$-connected graph and every homotopically non-trivial simple closed curve
intersects the graph at least three times.
It follows that the face boundaries of $S$
form a disk system, say $\calc$. Suppose now that $S$ is a subdivision of a graph $G$ and that
$S'$ is another $G$-subdivision obtained from $S$ by rerouting
or triad exchange.
Then $S$ uniquely determines an
embedding of $S'$ in $\Sigma$ (up to a homotopic shift) and the disk system
 induced in $S'$ by
$\calc$ consists of the face boundaries in $S'$.

\newthm{polyhedral} Let $G$ be an \afc\ graph, let $H$ be a graph, let $S$
be a $G$-subdivision in $H$, polyhedrally embedded in a surface
$\Sigma$, and assume that $S$ does not extend to an embedding of $H$.
Assume also that if $G$ has at most six vertices, then it is \ifc.
Then there exists a $G$-subdivision $S'$ in $H$ obtained from $S$
by repeated reroutings and at most one triad exchange
such that one of the following conditions holds for the induced
embedding of $S'$ into $\Sigma$:
\item{(i)} there exists an $S'$-path in $H$ such that no face
boundary of $S'$ includes both ends of the path,
\item{(ii)} there exists a free $S'$-cross on some face boundary
of $S'$, or
\item{(iii)} $H$ has an $S'$-separation, or
\item{(iv)} there exist an independent set $F\subseteq V(S')$ of size three,
a non-separating simple closed curve in $\Sigma$
intersecting $S'$ precisely in $F$, and an
$S'$-triad in $H$ with set of feet $F$ such that
$S'\backslash F$ is connected.
\smallskip

\proof Let ${\cal C}$ be the disk system described prior to the statement
of \refthm{polyhedral}.  By \refthm{summary2} there exists a $G$-subdivision
$S'$ in $H$ obtained from  $S$ by repeated reroutings and at most one
triad exchange that satisfies one of (i)--(v) of that lemma.
If (i), (ii) or (iii) holds, then our lemma holds.
Condition \refthm{summary2}(v) does not hold, because
$S$ does not extend to an embedding of $H$.
Thus we may assume that \refthm{summary2}(iv) holds.
Let $x_1,x_2,x_3$ be the feet of the
triad; since every pair of $x_1,x_2,x_3$ belong to a common face
boundary, there exists a simple closed curve $\phi$ passing through
$x_1,x_2,x_3$ and those faces.  Since no face boundary of $S'$
includes all of $x_1,x_2,x_3$ and
$S'\backslash \{x_1,x_2,x_3\}$ is connected, it follows that $\phi$ does not
separate $\Sigma$.  Thus (iv) holds.\qed

From now on we will be working exclusively with disk systems consisting
of peripheral cycles in subdivisions of $3$-connected planar graphs, and so
the notions such as $S$-jump or $S$-cross will refer to the disk system
consisting of all peripheral cycles.
If $G$ is planar, then there is no non-separating closed curve,
and hence condition (iv) from the above theorem cannot hold.
Thus we have the following corollary for planar graphs.
The corollary is used in~[\cite{DinOpoThoVerlarge}].

\newthm{polyplanar} Let $G$ be an \afc\ planar graph,
let $H$ be a non-planar graph, and let $S$
be a $G$-subdivision in $H$.
Assume also that if $G$ has at most six vertices, then $G$ is \ifc.
Then there exists a $G$-subdivision $S'$ in $H$ obtained from $S$
by repeated reroutings and at most one triad exchange
such that $S'$ and the disk system of peripheral cycles in $S'$
satisfy one of the following conditions:
\item{(i)} there exists an $S'$-path in $H$ such that no peripheral cycle
of $S'$ includes both ends of the path,
\item{(ii)} there exists a free $S'$-cross on some peripheral cycle
of $S'$, or
\item{(iii)} $H$ has an $S'$-separation.
\smallskip


Finally, we prove~\refthm{main}, which we restate in a slightly
stronger form.

\newthm{lemma3}Let $G$ be an \afc\ planar graph, let $H$ be an \afc\
non-planar graph, and let $S$ be a $G$-subdivision in $H$.
Assume that if $|V(G)|\le 6$, then $G$ is \ifc.
Then there exists a $G$-subdivision
$S'$ in $H$ obtained from $S$ by repeated reroutings and at most one
triad exchange such that $S'$ and the disk system of peripheral cycles in $S'$
satisfy one of the following conditions:
\item{(i)}there exists an $S'$-jump in $H$, or
\item{(ii)}there exists a free $S'$-cross in $H$ on some peripheral cycle
of $S'$.

\proof
Let $G,H,S$ be as stated. By~\refthm{polyplanar} there exists
a $G$-subdivision $S'$ in $H$ obtained from $S$
by repeated reroutings and at most one triad exchange such that
one of the conclusions of~\refthm{polyplanar} holds.
We may assume that $H$ has an $S'$-separation $(X,Y)$, for
otherwise the lemma holds.
Then $|X-Y|\ge2$, because $H[X]$ cannot be drawn in a disk with
$X\cap Y$ drawn on the boundary of the disk. The set $X-Y$ includes
at most one branch-vertex of $S'$ by the definition of $S'$-separation.
We claim that $|Y-X|\ge2$. This is clear if $S'$ has at least six branch-vertices;
otherwise $G$ has exactly five vertices and hence is \ifc.
It follows that $X$ cannot include four branch-vertices of $S'$, and so
$|Y-X|\ge2$, as claimed.
But that contradicts the \afcity\ of $H$.~\qed


We need the following lemma.
Let $G$ be a subdivision of a 3-connected planar graph, and let
$x,y$ be vertices or edges of $G$. We say that
$x,y$ are \dfn{cofacial} if some peripheral cycle in $G$ includes both $x$
and $y$.

\newthm{cofacial} Let $G$ be an \ifc\ planar graph, and let
$e\in E(G)$ and $v\in V(G)$ be not cofacial.  Then at least one end
of $e$ is not cofacial with $v$.

\proof Let us fix a planar drawing of $G$, and suppose for a
contradiction that both ends of $e$ are cofacial with $v$.
By~\refthm{periph} there exists a simple closed curve in the plane that
passes through $v$, the two ends of $e$, and is otherwise disjoint from $G$.
Since $v$ and $e$ are not cofacial this curve disconnects $G$,
contrary to the \ifcity\ of $G$.~\qed

We also need the following analogue of~\refthm{cofacial}.

\newthm{cofacial2} Let $G$ be an \ifc\ planar graph, and let
$e,f\in E(G)$  be not cofacial.  Then some  end
of $e$ is not cofacial with some end of $f$.

\proof
Let $u_1,u_2$ be the ends of $e$.
By~\refthm{cofacial} it suffices to show that one of $u_1,u_2$ is
not cofacial with $f$.
Thus we may assume for a contradiction that there exist
peripheral cycles $C_1,C_2$ in $G$, both containing $f$ and such that
$u_i\in V(C_i)$.
Let us fix a drawing of $G$ in the plane. By~\refthm{periph}
there exists a simple closed curve intersecting the graph $G$
three times: in $u_1$, $u_2$ and in an internal point of $f$.
However, that contradicts the \ifcity\ of $G$.~\qed

If we allow contracting edges and $G$ has no
peripheral cycles of length three,
then \refthm{lemma3} can be further simplified.
The next result is a restatement of~\refthm{maincor},
because every triangle-free \afc\ graph is \ifc.

\newthm{minor}Let $G$ be a triangle-free \ifc\ planar graph,
and let $H$ be an \afc\
non-planar graph such that $H$ has a subgraph isomorphic to a subdivision
of $G$. Then there exists a graph $G'$ such that $G'$ is isomorphic
to a minor of $H$, and either
\item{(i)}$G'=G+uv$ for some vertices $u,v\in V(G)$ such that
no peripheral cycle of $G$ contains both $u$ and $v$, or
\item{(ii)}$G'=G+u_1v_1+u_2v_2$ for some distinct vertices
$u_1,u_2,v_1,v_2\in V(G)$ such that $u_1,u_2,v_1,v_2$ appear on
some peripheral cycle of $G$ in the order listed.

\proof By \refthm{lemma3} there exist a \he\ $\eta:G\emb S\subseteq H$
and either an $S$-jump or a free $S$-cross.  Assume first that $P$ is
an $S$-jump with ends $a$ and $b$.  If both $a$ and $b$ are
branch-vertices, then (i) holds.  Let us assume that $a$ is a
branch-vertex, say $a=\eta (v)$ and that $b$ belongs to the interior of
$\eta (e)$ for some edge $e\in E(G)$.
Since $P$ is an $S$-jump it follows that $v$ and
$e$ are not cofacial.  By \refthm{cofacial} there exists an end
$u$ of $e$ such that $u$ and $v$ are not cofacial.  Then
$G+(u,v)$ satisfies (i).  We may therefore assume that neither
$a$ nor $b$ is a branch-vertex.  Let $a$ be an internal vertex of
$\eta (f)$ and let $b$ be an internal vertex of $\eta (e)$,
where $e,f\in E(G)$ are not cofacial.
By~\refthm{cofacial2} there is an end $u$ of $e$ that is not cofacial
with an end $v$ of $f$.
It follows that $G+(u,v)$ satisfies (i).

\junk{
as
before.  Let the ends of $f$ be $u_1$ and $u_2$, and let the ends of $e$
be $u_3$ and $u_4$.  If there exist integers $i\in \{1,2\}$ and
$j\in\{3,4\}$ such that $u_i$ and $u_j$ are not cofacial, then
$G+u_iu_j$ satisfies (i), and so we may assume that $u_i$ and
$u_j$ are cofacial for all $i=1,2$ and $j=3,4$.  By \refthm{cofacial} both
$u_1,u_2$ are cofacial with $e$ and both $u_3,u_4$ are cofacial
with $f$.  Thus there exist peripheral cycles $C_1,C_2,C_3,C_4$ such
that for $i=1,2,3,4$ the cycle $C_i$ includes $\{u_1,u_2,u_3,u_4\}
-\{u_i\}$.  By \refthm{periphcap} the graph $G$ is isomorphic to
$K_4$, contrary to \ifcity.  This completes the case where $H$ has
an $S$-jump.
} 

We may therefore assume that $P_1,P_2$ is a free $S$-cross in $H$ on some
peripheral cycle $\eta(C)$ of $S$, where $C$ is a peripheral cycle in $G$.  
Let $U$ be the set of feet  of this
cross, and let $B=V(C)$.
We define a bipartite graph $J$ with bipartition $(U,B)$ by
saying that $u\in U$ is adjacent to $b\in B$ if some subpath of $\eta(C)$ has
ends $u$ and $\eta(b)$ and includes no vertex of $U\cup \eta(B)$ in its interior.
Since $C$ has at least four vertices and
$P_1,P_2$ is a free cross, Hall's
theorem implies that $J$ has a complete matching from $U$ to $B$.
Let $U$ be matched into $\{u_1,u_2,v_1,v_2\}$,
where  $u_1,u_2,v_1,v_2$
occur on $C$ in the order listed.  Then $G+u_1v_1+u_2v_2$
satisfies (ii), as desired.~\qed

\newsection{appl}AN APPLICATION

By the \dfn{cube} we mean the graph of the $1$-skeleton of the $3$-dimensional
cube. As an application of the results of this paper
 we examine non-planar graphs
that have a subgraph isomorphic to a subdivision of the cube.
Other applications  appeared in [\cite{DinOpoThoVerlarge}, \cite{ThoThoTutte}].
Let $W$ denote the graph obtained from the cube by adding an edge joining
two vertices at distance three, and let $V_8$ be the graph obtained
from a cycle of length eight by adding  edges joining every pair of
diagonally opposite vertices. See Figure~\reffig{twographs}.

\goodbreak\midinsert
\bigskip
\bigskip
\centerline{\includegraphics[scale=0.7]{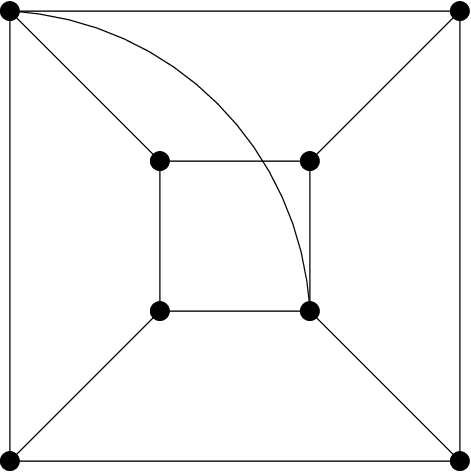}\quad\quad\includegraphics[scale=0.7]{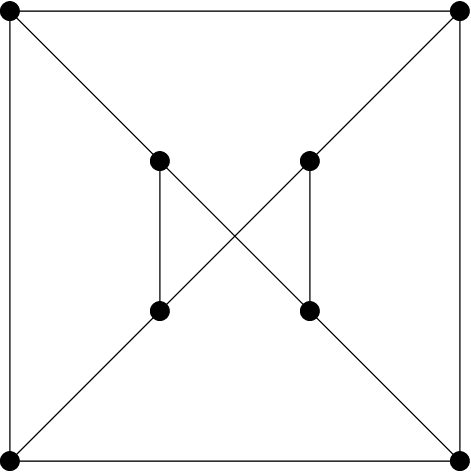}}
\bigskip
\centerline{\bf Figure~\newfig{twographs}.~~The graphs $W$ and $V_8$.}
\bigskip
\endinsert


\newthm{cube}Let $H$ be an \afc\ non-planar graph that has a subgraph isomorphic
to a subdivision of the cube. Then $H$ has a subgraph isomorphic to
a subdivision of $V_8$ or $W$.

\proof Let $K$ denote the cube.
By \refthm{lemma3} there exists a
\he\ $\eta:K\emb S\subseteq H$ such that (i) or (ii) of~\refthm{lemma3} holds.
Suppose first that (i) holds,
and let $P$ be a path as in (i) with ends $x$ and $y$. If $\eta(u)=x$
and $\eta(v)=y$, where $u,v\in V(K)$ are at distance three in $K$, then
$\eta$ can be extended to yield a $W$-subdivision in $H$, and the theorem holds.
Otherwise it is easy to see that $\eta$ can be extended to produce a
$V_8$-subdivision in $H$.

If (ii) holds, then there exists a free $\eta$-cross in some cycle $\eta(C)$
of $S$, where $C$  is a cycle of $K$ of length four.
Let the vertices of $C$ be $v_1,v_2,v_3,v_4$, in order.
Let $K'$ be obtained from $K$ by deleting the edges $v_1v_2$ and $v_3v_4$,
and adding the edges $v_1v_3$ and $v_2v_4$. The existence of the free
cross implies that $H$ has a subgraph isomorphic to a subdivision of
$K'$. But $K'$ is isomorphic to $V_8$, and so the result holds.~\qed

Theorem~\refthm{cube} is one step in the proof of the following
beautiful theorem of Maharry and Robertson~[\cite{MahRob}].

\newthm{V8} Let $G$ be an \ifc\ graph with no minor
isomorphic to $V_8$.  Then $G$ satisfies one of the following conditions:
\item{\rm (i)} $G$ has at most seven vertices,
\item{\rm (ii)} $G$ is planar,
\item{\rm (iii)} $G$ is isomorphic to the line graph of $K_{3,3}$,
\item{\rm (iv)} there is a set $X\subseteq V(G)$ of at most four
vertices such that $G\backslash X$ has no edges,
\item{\rm (v)} there exist two adjacent vertices $u,v\in V(G)$ such
that $G\backslash u\backslash v$ is a cycle.

Let $G$ be an \ifc\ graph on at least eight vertices.
In the first step Maharry and Robertson show that $G$
either is isomorphic to the line graph of $K_{3,3}$,
or has two disjoint cycles, each of length at least four.
Thus we may assume the latter, in which case \ifcity\ implies that $G$
has a minor isomorphic to $V_8$ or the cube.
By~\refthm{cube} we may assume that $G$ has a subgraph isomorphic to
a subdivision of $W$.
Let $u,v$ be the vertices of $G$ that correspond to the two vertices of $W$ of
degree four, let $X',Y'$ be the two color classes of the bipartite graph $W$,
and let $X$ and $Y$ be the sets of vertices of $G$ that correspond to 
$X'$ and $Y'$, respectively. 
Now it remains to show that either $G\backslash\{u,v\}$ is a cycle,
or that $G\backslash X$ or  $G\backslash Y$ has no edges.
To this end one can profitably apply the result of~[\cite{JohThoGener}].
We omit the details.

\newsection{improve}IMPROVING LEMMA~\refthm{planarlemma}

The objective of this section is to prove~\refthm{planarlemma2},
a version of~\refthm{planarlemma} that does not use rerouting.
Recall that since after~\refthm{polyhedral} all disk systems consist of peripheral cycles
of subdivisions of $3$-connected planar graphs.
The following is a version of~\refthm{localtriad} that uses no
rerouting or  triad exchange.

\newthm{triadjump} Let $G,H$ be graphs, where $G$ is \ifc\ and planar and
is not isomorphic to the cube,
let $S$ be a $G$-subdivision
in $H$, let $\calc$ be the disk system of peripheral cycles in $S$, and let
$Q_1,Q_2,Q_3$ be a local $S$-triad in $H$ centered at $v\in V(S)$
such that the $S$-bridge containing
$Q_1\cup Q_2\cup Q_3$ has an attachment $y$ that does not
belong to any of the segments incident with $v$.  Then there
exists an $S$-jump in $H$ with one end $y$.

\proof
Let $Z_1,Z_2,Z_3$ be the three segments incident with $v$, let $v_i$
be the other end of $Z_i$, and let $x_1,x_2,x_3$ be the feet of the
triad $Q_1,Q_2,Q_3$ numbered so that $v_i\in V(Z_i)$.
Let us fix a drawing of $S$ in the sphere.
For distinct integers $i,j,k\in\{1,2,3\}$ let $f_i$ be the face of $S$
incident with $Z_j$ and $Z_k$.
By hypothesis there exists a
path $P$ with ends $x\in V(Q_1\cup Q_2\cup Q_3)-\{x_1,x_2,x_3\}$
and $y\in V(S)-V(Z_1\cup Z_2\cup Z_3)$, disjoint from $S\backslash y$.
We may assume that for all $i=1,2,3$ the vertices $y$ and $x_i$ are
incident with the same face of $S$, for otherwise the lemma holds;
let $g_i$ denote that face.
Let $i,j,k\in\{1,2,3\}$ be distinct.
There exists a simple closed curve $\phi_k$ that intersects $S$ in
$\{y,v_i,v_j\}$
and is otherwise contained in $g_i\cup g_j\cup f_k$.
By the \ifcity\ of $G$ one of the disks bounded by $\phi_k$ includes
at most one branch-vertex of $S$.
Let $u_k$ denote that branch-vertex if it exists; otherwise $u_k$ is
undefined and $g_i=g_j=f_k$.
It follows that $G$ has at most eight vertices; the corresponding
branch-vertices of $S$ are $v,v_1,v_2,v_3,y$ and a subset of $\{u_1,u_2,u_3\}$.
Since $v_1$ has degree at least three, we deduce that at least one of
$u_1,u_2,u_3$ exists, say $u_3$ does. Then no segment of $S$ has ends
$y$ and $v_1$, or $y$ and $v_2$, or $v_1$ and $v_3$, or $v_2$ and $v_3$,
by the \ifcity\ of $G$.
Since $v_1$ and $v_2$ have degree at least three, it follows that
$u_1$ and $u_2$ also exist.
Thus $G$ is isomorphic to the cube, a contradiction.~\qed




We need to prove a variant of~\refthm{planarlemma}, where rerouting is not
used.
First we need \xx{two definitions}.
\xx{%
Let $S$ be a subdivision of a $3$-connected planar graph, let ${\cal C}$ be the disk system in $S$ consisting of peripheral cycles, let $Z_1,Z_2$ be distinct segments of $S$ with a common end $v$ such that they are both subgraphs of a disk $C\in\cal C$ and for $i=1,2$ let $v_i$ be the other end of $Z_i$. Let $P_1,P_2,P_3$ be paths such that
\item{$\bullet$} the ends of $P_i$ are $x_i$ and $y_i$,
\item{$\bullet$} $v_1,x_1,x_3,v,y_3,y_1,v_2$ appear on $Z_1\cup Z_2$ in the order listed, where possibly $v_1=x_1$ and/or $v_2=y_1$, but all other pairs are distinct,
\item{$\bullet$} $x_2$ is an internal vertex of $P_1$ and $y_2=v$,
\item{$\bullet$} the paths $P_1,P_2,P_3$ share no internal vertices with each other or with $S$, and 
\item{$\bullet$} the $S$-bridges containing $P_1$ and $P_3$ have no attachments outside $Z_1\cup Z_2$.}

\xx{%
\noindent
In those circumstances we say that $P_1,P_2,P_3$ is an {\sl $S$-leap}. }

Let $S$ be a subdivision of a $3$-connected planar graph,
let $W$ be a segment of $S$, let $z,w$ be the ends of $W$,
and let $P_1,P_2$ be two disjoint $S$-paths in $H$
with ends $x_1,y_1$ and $x_2,y_2$, respectively, such that
$z,x_1,x_2,y_1,w\in V(W)$ occur on $W$ in the order listed,
and $y_2\not\in V(W)$. Let $P_3$ be a path disjoint from
$V(S)-\{y_2\}$ with one end  $x_3\in V(P_1)$
and the other  $y_3\in V(P_2)$ and otherwise disjoint from
$P_1\cup P_2$.
Thus $P_1,P_2,P_3$ is an $S$-tripod based at $W$.
Let $\calc$ be the disk system in $S$ consisting of peripheral cycles, and
let $C,C'$ be the two disks that contain $W$.  Let
$y_2\in V(C)-V(C')$, and let $P_4$ be an $S$-path with
ends $x_4$ and $y_4$, where $x_4$ belongs to the interior of $x_1Wy_1$
and $y_4\in V(C')-V(C)$, such that no internal vertex of $P_4$
belongs to $P_1\cup P_2\cup P_3$.
For $i=1,2$ let $B_i$ be the $S$-bridge of $H$ that includes $P_i$.
Let us assume further that
\item{$\bullet$}
all attachments of $B_1$ and $B_2$ belong
to $C$,
\item{$\bullet$}
every $S$-bridge other than $B_1$
or $B_2$ that has an attachment in the interior of $x_1Wy_1$ has
all its attachments in $V(C')\cup\{y_2\}$, and
\item{$\bullet$}if $B_1\ne B_2$, then for $i=1,2$ the vertex $y_2$ is the only
attachment of $B_i$ that does not belong to $W$.

\noindent
In those circumstances we say
that the quadruple $P_1,P_2,P_3,P_4$ is an $S$-\dfn{tunnel}.
It is worth noting that if $B_1\ne B_2$, then $y_2=y_3$.
See Figure~\reffig{tunnel}.

\goodbreak\midinsert
\bigskip
\bigskip
\centerline{\includegraphics[scale=1.70]{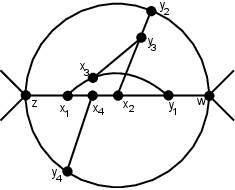}}
\bigskip
\centerline{Figure~\newfig{tunnel}. An $S$-tunnel.}
\bigskip
\endinsert

\newthm{planarlemma2}Let $G$ be an \ifc\ planar  graph,
let $H$ be
a graph, and let $S$ be a $G$-subdivision
in $H$
such that every unstable $S$-bridge is $2$-separated from $S$.
Then one of the following conditions holds:
\item{(i)}there exists an $S$-jump, or
\item{(ii)}there exists a weakly free $S$-cross in $H$, or
\item{(iii)}$H$ has an $S$-separation $(X,Y)$ such that $X-Y$ includes
no branch-vertex of $S$, or
\item{(iv)}there exists an $S$-triad, or
\item{(v)} there exists an $S$-tunnel, or
\item{(vi)} the graph $H$ is planar, \xx{or
\item{(vii)} there exists an $S$-leap.}
\smallskip

\proof
We proceed by induction on $|V(G)|$.
Suppose for a contradiction that
none of (i)--(vi) holds.
As in Claim (2) of~\refthm{planarlemma} we may assume that every $S$-bridge
is \stable, for otherwise the lemma follows by induction.
Now since every  $S$-bridge is \stable,
it follows from~\refthm{claim1}
and the fact that (i) and (iv) do not hold that
for every $S$-bridge $B$ there exists a unique disk $C$
such that all attachments of $B$ belong to $V(C)$.
 For every disk $C$ of $S$ let $H_C$ be the union
of $C$ and all $S$-bridges $B$ whose attachments are included
in $V(C)$. Since (vi) does not hold,
there exists a disk $C$ of $S$ such that $H_C$
does not have a planar drawing with $C$ bounding the infinite
face.
The same argument as in the proof of Claim~(3) of~\refthm{planarlemma}
shows that there exists \xx{either} an $S$-tripod \xx{or an $S$-leap.
(At the very end of the proof of Claim~(3) of~\refthm{planarlemma}
rerouting is used. Instead we get an $S$-leap.)}



Let us select a segment $Z$ and vertex $y_2\not\in V(Z)$ such that there
exists an $S$-tripod $P_1,P_2,P_3$ based at $Z$ with feet $x_1,y_1,x_2,y_2$.
Let $x_3\in V(P_1)$ and $y_3\in V(P_2)$ be the ends of $P_3$; we say
that $y_3P_2y_2$ is the \dfn{leg} of the $S$-tripod.  Let us, in
addition, select an $S$-tripod based at $Z$ so that its leg is minimal.
Let the leg be $L$.
We say that
a vertex $z\in Z$ is \dfn{sheltered} if $z$ is an internal vertex
of $x_1Zy_1$ for some $S$-tripod based at $Z$ with feet $x_1,y_1,x_2,y_2$
and leg $L$, and we say that the tripod \dfn{shelters} the vertex $z$.
Now let $x'_1,y'_1\in V(Z)$ be not sheltered but such that every
internal vertex of $x'_1Zy'_1$ is sheltered, and let $X'$ be the
union of $x'_1Zy'_1$ and $V(P_1\cup x_2P_2y_3\cup P_3)$, over all $S$-tripods
$P_1,P_2,P_3$ with leg $L$ that shelter an internal vertex of $x'_1Zy'_1$.

Let
$Y'=V(S\cup L)-(X'-\{x'_1,y'_1,y_3\})$. If there is no path between
$X'$ and $Y'$ in $H\backslash \{x'_1,y'_1,y_3\}$, then there exists
a separation $(X,Y)$ of order three in $H$ with $X'\subseteq X$
and $Y'\subseteq Y$ (and hence $X\cap Y=\{x'_1,y'_1,y_3\}$). Then
$(X,Y)$ is an $S$-separation, and hence (iii) holds, a contradiction.
Thus there exists a path
$P$ in $H\backslash \{x_1,y_1,y_3\}$ with ends $x\in X'$ and
$y\in Y'$.
We may assume that $P$ has no internal vertex in $X'\cup Y'$;
thus $P$ has no internal vertex in $S$.
If $x\in V(Z)$ let $P_1,P_2,P_3$ be an $S$-tripod that shelters
$x$; otherwise let $P_1,P_2,P_3$ be an $S$-tripod that shelters some
vertex of $x'_1Zy'_1$ such that $x\in V(P_1\cup P_2\cup P_3)$.

We claim that $P$ may be chosen so that $y\not\in V(Z)\cup V(L)$.
It is clear that $y\not\in V(L)$ by the choice of $L$, and so we
may assume that $y\in V(Z)$.  Let $B$ be the $S$-bridge that
includes $P$.  If $B$ includes at least one of $P_1,P_2,P_3$, then
$B\cup P_1\cup P_2\cup P_3$ includes an $S$-tripod that shelters
$x'_1$ or $y'_1$, a contradiction.  The same conclusion holds (or we obtain
contradiction to the minimality of $L$) if the only attachment of $B$
outside $Z$ is $y_2$.  Thus we may assume that $B$ has an attachment
in $V(S)-V(Z)-\{y_2\}$, and so $P$ may be replaced by a path with an end not in
$V(Z)\cup V(L)$.  This proves our claim that we may assume that
$y\not\in V(Z)\cup V(L)$.

Let $C_1,C_2$ be the two disks of $S$ that include $Z$.
Then $y_2\in V(C_1\cup C_2)$,
for otherwise $P_2$ is an $S$-jump and (i) holds.
From the symmetry we may assume that $y_2\in V(C_1)$.
Thus $y_2\not\in V(C_2)$ by (X1).
Assume first that $x\not\in V(Z)$.
If $y\in V(C_1)$, then $P_1\cup P_2\cup P_3\cup P$ includes a weakly free cross
on $C_1$, a contradiction.
Thus $y\not\in V(C_1)$.
Since for every $S$-bridge there is a disk that includes all the attachments
of the $S$-bridge, there exists a disk
$C_3\in\calc$ such that either
$x_2,y_2,y\in V(C_3)$ or $x_1,y_1,y_2,y\in V(C_3)$.
But $y\not\in V(C_1)$, and hence $C_3\ne C_1$. But $C_1,C_2$ are the
only two disks that contain $x_2$ and the only two disks that contain
both $x_1$ and $y_1$. Thus $C_3=C_2$, contrary to $y_2\not\in V(C_2)$,
a contradiction which completes the case $x\not\in V(Z)$.
We may therefore assume that $x\in V(Z)$, and that $P$ cannot be chosen with
$x\notin V(Z)$.  If $y\notin V(C_1)\cup V(C_2)$, then
$P$ is an $S$-jump, contrary to the fact that  (i) does not hold,
and if $y\in V(C_1)$, then $P\cup P_1\cup P_2
\cup P_3$ includes a weakly free $S$-cross,
contrary to the fact that (ii) does not hold.
Thus $y\in V(C_2)-V(C_1)$.
We claim that $P_1,P_2,P_3,P$ is an $S$-tunnel.
To this end let  $B_i$ be the $S$-bridge containing
$P_i$ for $i=1,2$.
The fact that the case $x\not\in Z$ and $y\not\in V(C_1)$ earlier in this
paragraph led to a contradiction implies that all attachments of $B_1$ and
$B_2$ belong to $C_1$.
Furthermore, if $B_1\ne B_2$ and one of them has an attachment in
$V(C_1)-V(Z)-\{y_2\}$, then $B_1\cup B_2$ includes a weakly-free cross,
contrary to the fact that (ii) does not hold.
Finally, let $B$ be an $S$-bridge other than $B_1$ or $B_2$ with
an attachment in the interior of $x_1Zy_1$.
The argument in this paragraph for the case $x\in V(Z)$ shows that
every attachment of $B$ belongs to $V(C_2)\cup\{y_2\}$.
This proves that $P_1,P_2,P_3,P$ is an $S$-tunnel, as desired.~\qed


\newsection{apex}APEX GRAPHS

Let $G$ be a graph.
By a \dfn{mold} for $G$ we mean a collection $Z=(Z_e: e\in F)$ of
(not necessarily disjoint)
sets, where $F\subseteq E(G)$ and each $Z_e$ is disjoint from $V(G)$.
Given a mold $Z$ for $G$ we define a new graph $L$ as follows.
We add the elements of $\bigcup_{e\in F} Z_e$ to $G$ as new vertices.
We subdivide each edge $e\in F$ exactly once,
 denoting the new vertex by $\hat e$.
Finally, for every $e\in F$ and every $z\in Z_e$
we add an edge between  $z$ and $\hat e$.
We say that $L$ is the graph \dfn{determined by $G$ and $Z$.}

Assume now that there exists a \he\ of $L$ into a graph $H$,
assume that $G$ is planar, but that the graph obtained from $H$ by deleting
the vertices that correspond to $\bigcup_{e\in F} Z_e$ is not.
Can the results obtained thus far be extended to this scenario?
We will study this question in this section, and we will find that
under some simplifying assumptions the answer is yes.
The main technical lemma is~\refthm{mainapex}, from which
we derive~\refthm{apexcor}, the main result of this section.
When $|\bigcup_{e\in F} Z_e|=1$ the main result has a simpler
form, stated as~\refthm{oneapex}.

Actually, we will not be given a \he\ of $L$ into $H$, but
some hybrid between a \he\ and a minor containment instead.
We now introduce this hybrid.
Let us recall that $\eta:G\emb S\subseteq H$ means that $S$ is a
$G$-subdivision in $H$ and $\eta$ maps vertices of $G$ to vertices of $S$
and edges of $G$ to the corresponding paths of $S$.
Let $\eta:G\emb S\subseteq H$.
Let $Z=(Z_e: e\in F)$ be a mold for $G$.
We say that $Z$ is a \dfn{mold for $G$ in $H$} if
$Z_e\subseteq V(H)$ for every $e\in F$.
By abusing notation slightly we will regard $Z$ as a graph with
vertex-set $\bigcup_{e\in F} Z_e$ and no edges.
Thus we can speak of $(S\cup Z)$-bridges.
By an \dfn{$S\cup Z$-link} we mean a subgraph $B$ of $H$ such that
either $B$ is isomorphic to $K_2$ and its vertices but not its edge
belong to $S\cup Z$, or $B$ consists of a connected subgraph $K$ of
$H\backslash V(S\cup Z)$ together with some edges from $K$ to
$S\cup Z$ and their ends.
Thus every $S\cup Z$-bridge is an $S\cup Z$-link, but not
the other way around.
We say that a mold $Z$ is \dfn{feasible} for $\eta$ if for every $e\in F$
and every $z\in Z_e$ there exists an $S\cup Z$-link   $B_{ez}$ of $H$ such
that the following conditions hold for all $e\in F$ and all $z\in Z_e$:
\item{(i)}$Z_e\subseteq V(H)-V(S)$,
\item{(ii)}$z\in V(B_{ez})$,
\item{(iii)}$V(B_{ez})\cap V(B_{e'z'})\subseteq V(S\cup Z)$ for
all distinct $e,e'\in F$ and all $z\in Z_e$ and $z'\in Z_{e'}$,
\item{(iv)}either some internal vertex
of $\eta(e)$ belongs to $B_{ez}$, or
both  ends of $\eta(e)$ belong to $B_{ez}$ 
and $B_{ez}=B_{ez'}$ for all $z'\in Z_e$.



\noindent
If the mold $Z$ is feasible and the graphs $B_{ez}$ are as above, then
we say that the collection $(B_{ez}:e\in F,z\in Z_e)$
of graphs is a \dfn{cast} for $Z$ and $\eta$ in $H$.
Thus feasibility is the promised hybrid between \he s and minors,
as the next lemma explains.

\newthm{moldminor}
Let $G,H$ be graphs, let $Z=(Z_e:e\in F)$ be a mold for $G$ in $H$,
and let $L$ be the graph determined by $G$ and $Z$.
If  $\eta:G\emb S\subseteq H\backslash V(Z)$ and
$Z$ is feasible for $\eta$, then $L$ is isomorphic to a minor of $H$.
Conversely, if $\eta_0:L\emb S_0\subseteq H$ satisfies $\eta_0(z)=z$
for every $z\in V(Z)$ and $\eta$ is the restriction of $\eta_0$ to $G$,
then $Z$ is feasible for $\eta$.

\proof
Let $\eta:G\emb S\subseteq H\backslash V(Z)$ and
let $Z$ be feasible for $\eta$.
Thus there exists a cast $\Gamma=(B_{ez}:e\in F,z\in Z_e)$ for $Z$ and
$\eta$ in $H$.
For $e\in F$ we define a connected graph $K_e$ as follows.
If there exists $z\in Z_e$ such that $B_{ez}$ has no attachment in the
interior of $\eta(e)$, then let 
$B=B_{ez'}$ for all $z'\in Z_e$
(this exists by the last axiom in the definition of a feasible mold),
and let $K_e:=B\backslash V(S\cup Z)$.
Otherwise let
$K_e$ be the union of the interior of $\eta(e)$ and
$B_{ez}\backslash V(S\cup Z)$ over all $z\in Z_e$, and all edges from
the latter sets to the interior of $\eta(e)$.
Then for distinct edges $e,e'\in F$ the graphs $K_e$ and $K_{e'}$ are disjoint.
By contracting all but one edge of each path $\eta(e)$ for every
$e\in E(G)-F$ we obtain an $L$-minor, where each $v\in V(G)$
is represented by $\eta(v)$, each $z\in V(Z)$ is represented by itself,
and for $e\in F$ the vertex $\hat e$ of $L$ is represented by $K_e$.
Thus $H$ has an $L$-minor.

Conversely, if $\eta_0:L\emb S_0\subseteq H$ satisfies $\eta_0(z)=z$
for every $z\in V(Z)$ and $\eta$ is the restriction of $\eta_0$ to $G$,
then a cast for $Z$ and $\eta$ in $H$ is constructed by letting
$B_{ez}:=\eta_0(z\hat e)$.~\qed

A cast $(B_{ez}:e\in F, z\in Z_e)$ is \dfn{united} if there exist distinct edges
$e,e'\in F$ and (not necessarily distinct) vertices $z\in Z_e$
and $z'\in Z_{e'}$ such that $B_{ez}$ and $B_{e'z'}$ are subgraphs of the
same $S\cup Z$-bridge.
A cast $(B_{ez}:e\in F, z\in Z_e)$ is \dfn{full} if each $B_{ez}$ is
an $S\cup Z$-bridge.

\newthm{unitefull} 
Let $G,S,H$ be graphs, let $\eta:G\emb S\subseteq H$, and let
$Z=(Z_e: e\in F)$ be a feasible mold for $G$ in $H$.
Then there exists a cast for $Z$ and $\eta$ in $H$ that is either
united or full.

\proof
Let $\Gamma=(B_{ez}:e\in F, z\in Z_e)$ be a cast for $Z$ and $\eta$ in $H$.
If $\Gamma$ is not united, then we may replace each $B_{ez}$ by the
$S\cup Z$-bridge it is contained in, thereby producing a full cast.~\qed
\medskip

Let $G,S,H,\eta$ and $Z$ be as above.
As in earlier sections of the paper we will produce $S$-jumps and $S$-crosses.
However, in order for them to be useful we need them to behave well with
respect to a cast. That leads to the following definitions.
An $S$-path $P$ is \dfn{compatible} with a full cast $(B_{ez}:e\in F,z\in Z_e)$
if $P$ is disjoint from $Z$ and
it is the case that if $P$ is a subgraph of $B_{ez}$ for some
$e\in F$ and $z\in Z_e$, then either
  one of the ends of $P$ belongs to the interior of $\eta(e)$,
  or $B_{ez}$ has no attachment in the interior of $\eta(e)$
  (in which case both ends of $\eta(e)$ are attachments of $B_{ez}$ by
  the last axiom in the definition of cast)
  and one end of $P$ is an end of $\eta(e)$.
We say that a cross $P_1,P_2$ is
\dfn{compatible} with a full cast $(B_{ez}:e\in F,z\in Z_e)$
if it satisfies the following conditions
\item{(C1)} both $S$-paths $P_1,P_2$ are compatible with the cast,
\item{(C2)} if $P_1,P_2$ are
subgraphs of the same $S$-bridge $B$, then $B=B_{ez}$ for no $e\in F$
and $z\in Z_e$,
\item{(C3)} there exists an index $i \in \{1,2\}$ such that $P_i$ has no attachments
in the interior of $\eta(e)$ for any $e \in F$.

\newthm{compatpath}
Let $G,H$ be graphs, let $\eta:G\emb S\subseteq H'$ be a \he,
let $Z=(Z_e:e\in F)$ be a feasible mold for $G$ in $H$, let
$(B_{ez}:e\in F,z\in Z_e)$ be a full cast for $Z$ and $\eta$ in $H$,
and let $P$ be an $S$-path in $H\backslash Z$.
Let $F'$ be the set of all edges $e\in F$ such that if
$P$ is a subgraph of $B_{ez}$ for some $z\in Z_e$, then either one end of $P$
is an internal vertex of the path $\eta(e)$, or
$B_{ez}$ has no attachment in the interior of $\eta(e)$
  and one end of $P$ is an end of $\eta(e)$.
Then $P$ is compatible with the cast $(B_{ez}:e\in F',z\in Z_e)$.
\smallskip

\noindent
The proof is clear.
\smallskip

The following lemma shows how to use $S$-jumps compatible
with a full cast.

\newthm{compatpathminor}
Let $G$ be an \ifc\ planar graph, let $H$ be a graph,
let $Z=(Z_e:e\in F)$ be a mold for $G$ in $H$,
and let $L$ be the graph determined by $G$ and $Z$.
If  $\eta:G\emb S\subseteq H\backslash V(Z)$ is a \he, $\Gamma$ is a full cast
for $Z$ and $\eta$ in $H$ and there exists an $S$-jump compatible
with $\Gamma$, then there exist vertices $u,v\in V(L)-V(Z)$
that are not cofacial in $L\backslash V(Z)$ such that
$L+uv$ is isomorphic to a minor of $H$.

\proof
Let $\Gamma=(B_{ez}:e\in F, z\in Z_e)$ and let $P$ be an $S$-jump
compatible with $\Gamma$.
The proof of~\refthm{moldminor} and the definition of compatible
path imply that there exists a graph $L'$
obtained from $L$ by subdividing at most two edges of $E(G)-F$
such that $L'+xy$ is isomorphic to a minor of $H$ for some two
vertices $x,y\in V(L')$ that are not cofacial in $L'\backslash V(Z)$.
This is straightforward, except for the case when $P$ is a subgraph
of $B_{ez}$ for some $e\in F$ and  $z\in Z_e$, and $B_{ez}$ has
no attachment in the interior of $\eta(e)$.
Then one end of $P$, say $\eta(x)$, is an end of $\eta(e)$ by the
definition of compatible path.
If the other end of $P$ is $\eta(y)$ for some $y\in V(G)$, then
$L+y\hat e$ is isomorphic to a minor of $H$, and $y$ and $\hat e$ are
not cofacial in $L\backslash V(Z)$, because $y$ and $x$ are not.
If the other end of $P$ belongs to the interior of $\eta(f)$ for some
$f\in E(G)$, then similarly
$L'+w\hat e$ is isomorphic to a minor of $H$, where $L'$ is obtained
from $L$ by subdividing the edge $f$ and $w$ denotes the new vertex,
and $w$ and $\hat e$ are not cofacial in $L'\backslash V(Z)$.
This completes the argument that $L'+xy$ is isomorphic to a minor of $H$.

The conclusion now follows from~\refthm{cofacial} and~\refthm{cofacial2}
in the same way as~\refthm{minor}.~\qed

Next we show how to use $S$-crosses compatible
with a full cast.

\newthm{compatcrossminor}
Let $G$ be an \ifc\ triangle-free planar graph, let $H$ be a graph,
let $Z=(Z_e:e\in F)$ be a mold for $G$ in $H$,
and let $L$ be the graph determined by $G$ and $Z$.
If  $\eta:G\emb S\subseteq H\backslash V(Z)$ is a \he, $\Gamma$ is a full cast
for $Z$ and $\eta$ in $H$ and there exists a free $S$-cross compatible
with $\Gamma$, then either
\item{(i)}there exist vertices $u,v\in V(L)-V(Z)$
that are not cofacial in $L - V(Z)$ such that
$L+uv$ is isomorphic to a minor of $H$, or
\item{(ii)}there exists a facial cycle $C$ of $L-V(Z)$
and distinct vertices $u_1,u_2,v_1,v_2\in V(C)$ appearing on $C$ in the
order listed such that $L+u_1v_1+u_2v_2$ is isomorphic to a minor of $H$, and for $i=1,2$
the vertices $u_i$ and $v_i$ are not adjacent in $G$.

\proof 
Let $\Gamma=(B_{ez}:e\in F, z\in Z_e)$,
let $P_1,P_2$ be a free $S$-cross compatible with $\Gamma$ and let $C$ 
be the facial cycle of $G$ such that $\eta(C)$ contains the feet of this cross.
We begin by considering the case when, for some $i \in \{1,2\}$, the path  $P_i$
is a subgraph of $B_{ez}$ for some $e\in F - E(C)$ and $z\in Z_e$.
Let  $P_i$  have ends $x,y\in V(\eta(C))$;
then $x,y$ do not belong to the same segment of $S$.
Condition (C1) in the definition of compatible path
implies that $B_{ez}$ has no attachment in the interior of $\eta(e)$
and one of $x,y$ is an end of $\eta(e)$, say $y$ is its end.
If $x$ is a branch vertex of $S$, then let $u\in V(G)$ be such that
$\eta(u)=x$; otherwise let $u\in V(G)$ be such that $y$ and $\eta(u)$ do not belong to the same segment. Such a choice is possible because $x$ and $y$ do not belong to the same segment. Lemma \refthm{periphcap} implies that
$u$ and $\hat e$ are not cofacial in $L\backslash V(Z)$.
The presence of $P_i$ guarantees that $L+u\hat e$ is isomorphic to
a minor of $H$, and so (i) holds.

Thus we may assume that if $P_i$ is a subgraph of $B_{ez}$ for some $e\in F$ and $z\in Z_e$, then
$e \in E(C)$. As the cross is free and condition (C3) in the definition of a cross compatible with a cast is satisfied, at most one foot
of the cross $P_1,P_2$ belongs to the interior of $\eta(e)$ for every $e \in F$.  

\junk{
We may now assume that there exists a graph $L'$
obtained from $L$ by subdividing edges of $C'$ and a free cross in the graph $L'-V(Z)$ with feet on the cycle $C$ obtained from $C'$ via this subdivision. As in~\refthm{compatpathminor},
this is easy unless, $P_1$ or $P_2$ is a subgraph of  $B_{ez}$ for some $e\in F$ and $z\in Z_e$, such that  $B_{ez}$ has no attachment in the interior of $\eta(e)$. Suppose $P_1$ is. Then $P_2$ is not a subgraph of $B_{ez}$ by (C2).  If  $P_2$ has an attachment $x$ in the interior of $\eta(e)$ and $y \not \in V(P_1)$ is an end of $\eta(e)$, then let $P_2'=P_2 \cup x\eta(e)y$. Let $P_2'=P_2$, otherwise. Contracting $B_{ez}$  to a single vertex and deleting the interior of $\eta(e) - P_2'$,  allows one to construct the graph $L'$ and the cross, as required, in this case. (See Figure~\reffig{cross95}.)
\vskip 10pt
\centerline{\includegraphics[scale=1.8]{cross95.eps}}
\bigskip
\centerline{Figure~\newfig{cross95}. Constructing a free cross in $L'-V(Z)$.}
\bigskip
}

We now repeat the argument of~\refthm{minor}, with slight modifications, and 
we also use the proof of~\refthm{moldminor}. 
We say that an edge $e\in F$ is \dfn{internal} if $B_{ez}$ has an attachment in the interior of $\eta(e)$
for every $z\in Z_e$, and otherwise we say that $e$ is \dfn{external}.
Let $U$ be the set of feet of the cross $P_1,P_2$, and let $B\subseteq V(L)$ consist of all vertices of $C$
and all vertices of the form $\hat e$, where $e\in F\cap E(C)$ is internal.
(Let us recall that $\hat e$ is the new vertex of $L$ that results from subdividing the edge $e$.)
We define a bipartite graph $J$ with bipartition $(U,B)$ as follows.
Let $u\in U$, and let us assume first that $u$ is a branch-vertex of $S$.
Let $i\in\{1,2\}$ be such that $u$ is a foot of $P_i$.
If $P_i$ is a subgraph of $B_{ez}$ for some external edge $e\in F$ and $z\in Z_e$ and 
$u$ is an end of $\eta(e)$, then we declare $u$ adjacent to $\hat e$ only.
Otherwise let $x\in V(C)$ be such that $u=\eta(x)$, and we declare
that $u$ is adjacent to $x$ only.
Thus we may assume that $u$ belongs to the interior of $\eta(e)$ for some $e\in E(C)$.
If $e\in F$ is an interior edge, then we declare $u$ to be adjacent to $\hat e$ only.
Otherwise $u$ will be adjacent to every end $x$ of $e$ such that the subpath $Q$ of $\eta(e)$
between $\eta(x)$ and $u$ includes no member of $U$ in its interior.
In that case we say that $Q$ \dfn{represents} the edge $ux$ of $J$.
It follows similarly as in~\refthm{minor} that
the graph $J$ has a complete matching $M$ from $U$ to $B$,
but extra care is needed. In particular, we need condition (C2).
Furthermore, the matching $M$ may be chosen so that if $e=xy\in F\cap E(C)$ is internal,
then at least one of the vertices $x,\hat e,y$ is not saturated by $M$.

Let $U$ be matched by $M$ to the set $u_1,u_2,v_1,v_2 \in V(L)$, where 
$u_1,u_2,v_1,v_2$ appear on the cycle of $L$ that corresponds to $\eta(C)$ in the order listed.
We claim that $L+u_1v_1+u_2v_2$ is isomorphic to a minor of $H$.
Indeed, this follows similarly as in~\refthm{minor}, using the argument of the 
proof of~\refthm{moldminor}.
More specifically, we define the graphs $K_e$ as in the proof of~\refthm{moldminor}.
The proof of~\refthm{moldminor} shows that $L$ is isomorphic to a minor of $H$.
To obtain the same conclusion for  $L+u_1v_1+u_2v_2$ we make sure that
when contracting the edges of the paths $\eta(e)$  for $e\not\in F$ we contract
all edges of every  subpath of $\eta(e)$ that represents an edge of $M$.
We also need to contract all edges of paths that represent edges of $M$ and
are subpaths of $\eta(e)$ for external edges $e\in F$.
The path $P_i$ then gives rise to the edge $u_iv_i$.
%
If $e=u_iv_i \in E(G)$  for some $i \in \{1,2\}$,
then (since $u_1,u_2,v_1,v_2$ appear in the order listed)
$e\in F\cap E(C)$ is internal and one of $u_{3-i},v_{3-i}$ is equal to $\hat e$,
contrary to the choice of $M$.
Thus $u_i$ and $v_i$ are not adjacent in $G$ and (ii) holds.~\qed

Let $\eta:G\emb S\subseteq H$, and let $\eta':G\emb S'\subseteq H$
be obtained from $\eta$ by a rerouting.
If $H'$ is a subgraph of $H$ and both $S$ and $S'$ are subgraphs of $H'$,
then we say that the rerouting is \dfn{within} $H'$.

Our next objective is to give a sufficient condition for a rerouting to
preserve feasibility.
To that end we need to discuss the effect of reroutings on casts.
Let $G,H$ be graphs, let $\eta:G\emb S\subseteq H$ be a \he,
let $F\subseteq E(G)$,
and let $\eta':G\emb S'\subseteq H$ be obtained from $\eta$ by
a rerouting.
We say that the rerouting is \dfn{$F$-safe} if the following conditions
are satisfied:
\item{(i)}if the rerouting replaces a subpath of $S$ by an $S$-path $Q$
and $Q$ is a subgraph of an $S$-bridge $B$, and $e\in F$ is such that
either $B$ has an attachment in the interior of $\eta(e)$,
or both ends of $\eta(e)$ are attachments of $B$, then
the rerouting is an I-rerouting based at $\eta(e)$,
\item{(ii)}if the rerouting is a T-rerouting centered at $\eta(v)\in V(S)$,
then no edge of $G$ incident with $v$ belongs to $F$, and
\item{(iii)}if the rerouting is a V- or X-rerouting based at $\eta(e_1)$
and $\eta(e_2)$, then $e_1,e_2\not\in F$.

\noindent
Thus every proper I-rerouting is $F$-safe.

\newthm{reroutefeas}
Let $G,H$ be graphs, let $\eta:G\emb S\subseteq H$ be a \he,
let $Z$ be a mold for $G$ in $H$ that is feasible for $\eta$,
let there be a full cast for $Z$ and $\eta$ in $H$,
and let $\eta':G\emb S'\subseteq H$ be obtained from $\eta$ by an
$F$-safe rerouting within $H\backslash V(Z)$.
Then $Z$ is feasible for $\eta'$.

\proof
Let $Z=(Z_e:e\in F)$ and let $\Gamma=(B_{ez}:e\in F, z\in Z_e)$
be a full cast for $Z$ and $\eta$ in $H$.
Let $e\in F$ and $z\in Z_e$. We wish to define an $S'\cup Z$-link $B'_{ez}$.
If the rerouting is an I-rerouting, then let $W$ be a segment of $S$
such that the rerouting is based at $W$; otherwise let $W$ be the null graph.
The construction will be such that
$V(B'_{ez})\subseteq V(B_{ez}\cup W)$.
That will guarantee that the links thus defined will satisfy the
third axiom in the definition of feasibility.

Assume first that $B_{ez}$ includes an $S$-path $Q$ that replaced
a subpath $P$ of $S$ during the rerouting.
Since $\Gamma$ is a full cast, the $S$-bridge $B_{ez}$ either has an
attachment in the interior of $\eta(e)$,
or both ends of $\eta(e)$ are attachments of $B_{ez}$.
The first axiom in the definition of $F$-safety implies that the
rerouting is an  I-rerouting based at $\eta(e)$.
The $S\cup Z$-bridge $B_{ez}$ includes a path from $z$ to the interior of
$Q$; let $B'_{ez}$ be such a path with no internal vertex in $S'\cup Z$.
This completes the construction when $B_{ez}$ includes an $S$-path
that replaced a subpath $P$ of $S$ during the rerouting.

Thus we may assume that $B_{ez}$ includes no such $S$-path.
If no attachment of $B_{ez}$ belongs to $\eta(e)$
and to the interior of a subpath
of $S$ that got replaced by an $S$-path  during the rerouting,
then we let $B'_{ez}:=B_{ez}$.
We may therefore assume that an attachment $x$ of $B_{ez}$ belongs to
$\eta(e)$ and to the interior of a subpath
$P$ of $S$ that got replaced by an $S$-path $Q$ during the rerouting.
The second and third axiom in the definition of safety imply that
the rerouting is an I-rerouting and that $x$ belongs to the interior of
$\eta(e)$. Thus the I-rerouting is based at $\eta(e)$.
If the ends of $Q$ are not equal to the ends of
$\eta(e)$, then we define $B'_{ez}:=P\cup B_{ez}$.
It follows that the $S'\cup Z$-link $B'_{ez}$  has an attachment in the
interior of $\eta'(e)$.
Thus we may assume that $Q$ and $\eta(e)$ have the same ends.
In that case we define $B'_{ez}:=P\cup\bigcup_{z'\in Z_e}B_{ez'}$,
in which case $B'_{ez}=B'_{ez'}$ for all $z'\in Z_e$ and both ends
of $\eta'(e)$ are attachments of $B'_{ez}$.
Hence the $S'\cup Z$-links $B'_{ez}$ satisfy the last feasibility axiom.
The third axiom follows as indicated earlier, and other axioms are clear.

Thus $(B'_{ez}:e\in F,z\in Z_e)$ is a cast for $Z$ and $\eta'$ in $H$,
as required.~\qed

\newthm{properIreroute}
Let $G,H$ be graphs, let $Z=(Z_e:e\in F)$ be a mold for $G$ in $H$,
and let $\eta:G\emb S\subseteq H$ be a \he.
If $Z$ is feasible for $\eta$ and there exists a full cast for $Z$ and $\eta$
in $H$, and $\eta':G\emb S'\subseteq H$
is obtained from $\eta$ by a proper I-rerouting within $H\backslash V(Z)$,
then $Z$ is feasible for $\eta'$.

\proof
This follows immediately from~\refthm{reroutefeas}, because a proper
I-rerouting is $F$-safe.~\qed

The following is the main technical lemma of this section.

\newthm{mainapex}
Let $G$ be an \ifc\ planar graph not isomorphic to the cube,
let $H$ be a graph, and let $Z=(Z_e:e\in F)$ be a mold for $G$ in $H$.
Let $H':=H\backslash\bigcup_{e\in F} Z_e$, and let
$\eta_0:G\emb S_0\subseteq H'$ be a \he\ such that the mold $Z$ is feasible
for $\eta_0$.
Then there exist a \he\ $\eta:G\emb S\subseteq H'$ obtained from $\eta_0$
by repeated reroutings within $H'$ and a set $F'\subseteq F$ such that
every two edges in $F-F'$ are cofacial in $G$ and letting
$Z'$ denote the mold $(Z_e:e\in F')$
one of the following conditions holds:
\item{(i)}there is a united cast for $Z'$ and $\eta$, or
\item{(ii)}there exists an $S$-jump compatible with some full cast for
$Z'$ and $\eta$, or
\item{(iii)}there exists a free $S$-cross compatible with some full cast for
$Z'$ and $\eta$, or
\item{(iv)}$H'$ has an $S$-separation, or
\item{(v)}$H'$ is planar.

\proof
Let $\eta_0:G\emb S_0\subseteq H'$  and $Z$ be as stated.
We may assume that (i) does not hold.
We start with the following claim.

\claim{1} Let $\eta:G\emb S\subseteq H'$ be obtained from $\eta_0$
  by repeated reroutings within $H'$, let $F'\subseteq F$ be such that
  every two edges in $F-F'$ are cofacial in $G$, let
  $Z':=(Z_e:e\in F')$, and let there exist
  a full cast for $Z'$ and $\eta$ in $H'$.
  If $\eta'$ is obtained from $\eta$ by an $F'$-safe rerouting,
  then there is a full cast for $Z'$ and $\eta'$ in $H'$.

\noindent
To prove (1) we first notice that~\refthm{reroutefeas} implies that
$Z'$ is feasible for $\eta'$ in $H'$.
By~\refthm{unitefull} there is a cast for $Z'$ and $\eta'$ in $H'$
that is united or full.
The former does not hold by our assumption that (i) does not hold,
and hence the latter holds.
This proves (1).
\medskip

By~\refthm{stable} applied to the graphs $G,S_0$ and $H'$ there exists
a \he\ $\eta:G\emb S\subseteq H'$ obtained from $\eta_0$ by
repeated proper I-reroutings such that every unstable $S$-bridge
is $2$-separated from $S$.
Since every proper I-rerouting is $F$-safe,
it follows from (1) that there is a full cast for
$Z$ and $\eta$ in $H'$.
Let $\Gamma:=(B_{ez}:e\in F,z\in Z_e)$ be such a full cast.

\claim{2} If there exists an $S$-jump in $H'$, then the theorem holds.

\noindent
To prove (2) let $P$ be an $S$-jump.
If the $S$-bridge containing $P$ is
equal to $B_{ez}$ for some $e\in F$ and $z\in Z_e$,
then there is exactly one such edge $e$,
and we define $F':=F-\{e\}$; otherwise we let $F':=F$.
Then $P$ is compatible with $(Z_e:e\in F')$, and hence
$F'$ and $\eta$ satisfy (ii).
This proves (2).
\medskip


\claim{3} Let $u$ be a vertex of $G$ of degree three,
  let $F'$ be obtained from $F$ by removing all edges incident with $u$,
  let $\eta':G\emb S'\subseteq H'$ be obtained from $\eta$ by a sequence of
  $F'$-safe reroutings, and let there exist a local $S'$-triad centered
  at $\eta'(u)$.
  Then $F'$ satisfies the conclusion of the theorem.

\noindent
To prove (3) we first deduce from (1) that there exists a full
cast $\Gamma'=(B'_{ez};e\in F', z\in Z_e)$ for $Z'$ and $\eta'$.
Let $Z_1,Z_2,Z_3$ be the three segments of $S'$ incident
with $v:=\eta'(u)$, let $v_i$ be the other end of $Z_i$,
and let the local $S'$-triad be $Q_1,Q_2,Q_3$,
where $Q_i$ has end  $x_i\in V(Z_i)$.
Let $L_i:=v_iZ_ix_i$ and $P_i:=vZ_ix_i$.
We may assume that $\eta'$ and the triad $Q_1,Q_2,Q_3$ are chosen so that
$|V(L_1)|+|V(L_2)|+|V(L_3)|$ is minimum.




Let $X_1=V(P_1\cup P_2\cup P_3\cup Q_1\cup Q_2\cup Q_3)$ and
$Y_1=V(S)-(X_1-\{x_1,x_2,x_3\})$.
If $H'\backslash\{x_1,x_2,x_3\}$ has
no path between $X_1$
and $Y_1$, then
$H'$ has a separation $(X,Y)$ such that $X\cap Y=\{x_1,x_2,x_3\}$,
$X_1\subseteq X$, and $Y_1 \subseteq Y$.
Then $(X,Y)$ satisfies outcome (iv) of the theorem.

We may therefore assume that there exists a path $P$ in $H'$ as above. Let
the ends of $P$ be $x\in X_1- \{x_1,x_2,x_3\}$
and $y\in Y_1-\{x_1,x_2,x_3\}$.
We may assume that $P$ has no internal vertex in $X_1\cup Y_1$.
If $P$ is a subgraph of the $S'$-bridge
$B'_{ez}$ for some $e\in F'$ and $z\in Z_e$, then
we may assume that $y$ satisfies the following specifications.
If $B'_{ez}$ has an attachment in the interior of $\eta'(e)$, then
we may assume that $y$ belongs to the interior of $\eta'(e)$;
otherwise we may assume that $y$ is an end of $\eta'(e)$
(because both ends of $\eta'(e)$ are attachments of $B'_{ez}$ by the
last axiom in the definition of cast, and at least one end of $\eta'(e)$
does not belong to $Z_1\cup Z_2\cup Z_3$, because $G$ is \ifc).

Assume first that $x\in V(Q_1\cup Q_2\cup Q_3)$.
Then $y\notin V(L_1\cup L_2\cup L_3)$ by the choice of $Q_1,Q_2,Q_3$.
Since $G$ is not isomorphic to a cube we deduce from~\refthm{triadjump}
that there is an $S'$-jump with one end $y$.
The choice of $y$ implies that the $S'$-jump is compatible with $\Gamma'$,
and hence outcome (ii) holds.
This completes the case that $x\in V(Q_1\cup Q_2\cup Q_3)$.
Furthermore, it implies that we may assume that the $S'$-bridge
containing $Q_1\cup Q_2\cup Q_3$ has all attachments in $Z_1\cup Z_2\cup Z_3$.

Thus $x\in V(P_1\cup P_2\cup P_3)$.
%
Let $B$ be the
$S'$-bridge containing $P$.  If $B$ has an attachment outside
$Z_1\cup Z_2\cup Z_3$, then $P$ may be replaced by a path
with an end not in $Z_1\cup Z_2\cup Z_3$; otherwise replacing a path
of $P_1\cup P_2\cup P_3$ by $P$ is a T-rerouting centered at $v$, and
it is $F'$-safe.
The resulting \he\ has a triad that contradicts the choice of
$\eta'$ and $Q_1,Q_2,Q_3$.  Thus we may assume that $y\notin V(Z_1\cup
Z_2\cup Z_3)$.

We may assume that  $P$ is not an $S'$-jump, for
otherwise (ii) holds, because $P$ is compatible with $\Gamma'$ by
the choice of $y$.
Thus  there exists a disk $C$ in $S'$ such that $x,y\in V(C)$.
It follows that
$C$ includes two of the segments incident with $v$, say $Z_1$ and $Z_2$.
Now both $S'$-paths $Q_1\cup Q_2$ and $P$ are compatible with $\Gamma$,
the former because for $e\in F'$ the $S'$-bridge $B'_{ez}$ does not
include $Q_1\cup Q_2$, which in turn follows from
the fact that the $S'$-bridge
containing $Q_1\cup Q_2\cup Q_3$ has all attachments in $Z_1\cup Z_2\cup Z_3$.
Thus $Q_1\cup Q_2, P$ is an $S'$-cross compatible with $\Gamma$. (Condition (C3) holds as $Q_1 \cup Q_2$ does not have ends in the interior of images of edges in $F'$.)
The cross is free by the \ifcity\ of $G$.
This proves (3).

\claim{4} If there exists an $S$-triad in $H'$, then the theorem holds.

\noindent
To prove (4) assume that there exists an $S$-triad in $H'$.
The triad is local by~\refthm{planartriad}, and hence the claim follows
from (3) applied to the \he\ $\eta$.
This proves~(4).
\medskip

\xx{%
In preparation for the proof of (6) we prove the following special case.}

\claim{5}%
\xx{%
Assume that  there exist segments $Z_1,Z_2$ in $S$ with common end $v$ of degree 
at least four 
and the other ends $v_1,v_2$, respectively,
such that both are subgraphs of a disk $C$ and that there exists a weakly free $S$-cross $P_1,P_2$
 such that the ends of $P_i$ can
be labeled $x_i,y_i$ in such a way that $v_1,x_1,x_2,v,y_1,y_2,v_2$
occur on $Z_1\cup Z_2$ in the order listed.
Let $F'$ be obtained from $F$ by deleting all edges $e\in F$ such that
$\eta(e)$ is a subgraph of $C$,
and let $Z':=(Z_e:e\in F')$.
Then  $F'$ satisfies the conclusion of the theorem.
}

\noindent
\xx{To prove (5)}
let us define the \dfn{height} of the cross $P_1,P_2$ to be
$|E(L_1)|+|E(L_2)|$, where
$L_1:=v_1Z_1x_1$ and $L_2:=v_2Z_2y_2$.
We proceed similarly as in the proof of~\refthm{weakfree2}, but with
extra care.
Let $e_1,e_2\in E(G)$ be such that $\eta(e_i)=Z_i$.
Let $\eta_1:G\emb S_1\subseteq H'$ be a \he\ obtained from $\eta$ by
a sequence of proper V-reroutings based at $\eta(e_1),\eta(e_2)$ and let
$Q_1,Q_2$ be a weakly free $S_1$-cross based at $\eta_1(e_1),\eta_1(e_2)$
such that among all such triples $(\eta_1,Q_1,Q_2)$ this one minimizes
the height of the cross $Q_1,Q_2$.
Since every proper V-rerouting is $F'$-safe,
it follows from (1) that there is full cast for $Z'$ and $\eta_1$ in $H'$.
In order to prevent the introduction of unnecessary notation we now
make the assumption that $\eta=\eta_1$, $P_1=Q_1$ and $P_2=Q_2$.
This can be done with the proviso that for the remainder of the proof
of (5) $\Gamma$ is a full cast for $Z'$ and $\eta$ (as opposed to
a full cast for $Z$).


Let $X'$ be the vertex-set of
$P_1\cup P_2\cup vZ_1x_1\cup vZ_2y_2$
and let $Y'=V(S)-(X'-\{v,x_1,y_2\})$. If there is no path
in $H'\backslash\{v,x_1,y_2\}$ with one end in $X'$ and the other
in $Y'$, then there exists a separation $(X,Y)$ of order three with
$X'\subseteq X$ and $Y'\subseteq Y$. This separation satisfies (iv),
as required, and so we may assume that there exists a path $P$ in
$H'\backslash\{v,x_1,y_2\}$
with one end $x\in X'$ and the other end $y\in Y'$.

We first complete the proof of (5) assuming that $y\not\in V(Z_1\cup Z_2)$,
that at least one of $x,y$ is not in $V(C_1)$, and that at least one
of $x,y$ is not in $V(C_2)$.
From the symmetry we may assume that $x\in V(P_1\cup y_2Z_2v)$.
If $P_1$ is a subgraph of $B_{ez}$ for some $e\in F'$ and $z\in Z_e$,
then we may assume that either $y$ belongs to the interior of $\eta(e)$,
or $y$ is an end of $\eta(e)$ and $y\not\in V(C)$.
(This is indeed possible---by the choice of $F'$ at least one end of
$\eta(e)$ does not belong to $Z_1\cup Z_2$.)
If $y\not\in V(C\cup C_2)$, then $P_1\cup P$ includes an $S$-jump
with ends $y_1$ and $y$, which is compatible with $\Gamma$.
Thus (ii) holds.
Next let us assume that $y\in V(C_2)$. Then $y\not\in V(C)$, because
$C\cap C_2=Z_2$.
Since  $v$ has degree at least four,
(X2) implies that $V(C_1)\cap V(C_2)=\{v\}$.
It follows that $y\not\in V(C_1)$, and so
$P_1\cup P$ includes an $S$-jump, which is compatible with $\Gamma$.
Thus, again, (ii) holds.
We may therefore assume that $y\in V(C)$.
That implies that $P_1$ is a subgraph of $B_{ez}$ for no
$e\in F'$ and $z\in Z_e$, and so from the symmetry we may assume
the same about $P_2$.
Since $y\in V(C)$ we deduce that
$P_1\cup P_2\cup P$ includes a free cross.
Since $P_1,P_2$ are not subgraphs of any $B_{ez}$ for $e\in F'$ it follows that
the cross is compatible with $\Gamma$, and so (iii) holds.
This completes the case that $y\not\in V(Z_1\cup Z_2)$,
at least one of $x,y$ is not in $V(C_1)$, and at least one
of $x,y$ is not in $V(C_2)$.
Thus we may assume that the $S$-bridge that contains $P_1$
has all its attachments in $Z_1\cup Z_2$, and from the symmetry we
may assume the same about the $S$-bridge containing $P_2$.
In particular, the X-rerouting of $Z_1,Z_2$ that makes use of $P_1,P_2$
is proper, and hence is $F'$-safe.

Next we handle the case that $y\in V(Z_1\cup Z_2)$.
Again, from the symmetry we may assume that $x\in V(P_1\cup y_2Z_2v)$.
If $y\in V(L_1)$, then replacing $P_1$ by $P$ if $x\not\in V(P_1)$
and by $P\cup xP_1y_1$ otherwise produces a cross of smaller height,
contrary to
the choice of the triple $(\eta_1,Q_1,Q_2)$.
If $y\in V(L_2)$, then replacing $yZ_2x$ by $P$ if $x\not\in V(P_1)$
and replacing $yZ_2y_1$ by $P\cup xP_1y_1$ results in a \he\ $\eta'$
obtained from $\eta$ by a proper V-rerouting, and $P_1,P_2$ can be modified
to give a cross $P_1',P_2'$ such that the triple $(\eta',P_1',P_2')$
contradicts the choice of $(\eta_1,Q_1,Q_2)$. Thus $y\not\in V(Z_1\cup Z_2)$.

Finally, from the symmetry we may assume that $x,y\in V(C_2)$.
Let $B$ be the $S$-bridge containing $P$.
Since we may assume that $x,y$ cannot be chosen to satisfy any of
the cases already handled, it follows that every attachment of $B$
belongs to $C_2$.
Thus we may assume that if $B=B_{ez}$ for some $e\in F'$ and $z\in Z_e$,
then either $y$ is an internal vertex of $\eta(e)$, or
$B$ has no attachment in the interior of $\eta(e)$ and
$y$ is an end of $\eta(e)$.
Since
$V(C_1)\cap V(C_2)=\{v\}$, it follows that $y\not\in V(C_1)$.
Now let $\eta':G\emb S'\subseteq H'$ be obtained from $\eta$
by the X-rerouting using the cross
$P_1,P_2$, and let $Z_1',Z_2'$ be the segments of $S'$ corresponding
to $Z_1,Z_2$, respectively.
Thus $Z_1'=v_1Z_1x_1\cup P_1\cup y_1Z_2v$ and
$Z_2'=v_2Z_2y_2\cup P_2\cup x_2Z_1v$. (See Figure~\reffig{cross98}.)
As pointed out earlier, this rerouting is $F'$-safe.
It follows that $\Gamma$ is a cast for $Z'$ and $\eta'$.
Now $P$ is an $S'$-jump, and is compatible with $\Gamma$ by the
choice of $y$. Thus (ii) holds.
This completes the proof of (5).
\medskip

\centerline{\includegraphics[scale=1.5]{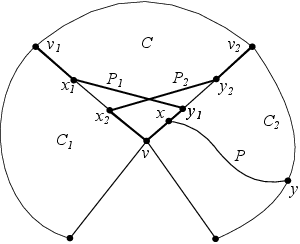}}
\bigskip
\centerline{Figure~\newfig{cross98}. $X$-rerouting in the proof of (5).}
\bigskip

\claim{\xx{6}} If there exists a weakly free $S$-cross in $H'$, then the theorem holds.

\noindent
To prove (\xx{6})
let $P_1,P_2$ be a weakly-free $S$-cross in $H'$ on a disk $C$, and assume for
a moment that the cross is free.
Let $F'$ be obtained from $F$ by removing the 
edges $e\in F$ such that $\eta(e) \subseteq C$.
The set $F'$ satisfies outcome
(iii) of the present theorem, unless, say, $P_1$ is
a subgraph of $B_{ez}$ for some $e\in F$ with $\eta(e)\not\subseteq E(C)$ and $z\in Z_e$.
But then the bridge $B_{ez}$ includes an $S$-jump or an $S$-triad in $H'$,
and hence the theorem holds by (2) and (4).
This concludes the case when $P_1,P_2$ is a free cross, and so we
may assume that it is not.

Thus there exist segments $Z_1,Z_2$ in $S$ with common end $v$
and the other ends $v_1,v_2$, respectively,
such that both are subgraphs of $C$ and such that the ends of $P_i$ can
be labeled $x_i,y_i$ in such a way that $v_1,x_1,x_2,v,y_1,y_2,v_2$
occur on $Z_1\cup Z_2$ in the order listed.
There are two cases depending on the degree of $v$.
\xx{If the degree of $v$ is at least four, then the claim follows from (5).
We may therefore assume}
that the degree of $v$ is three.
If the $S$-bridge of $H'$ that includes the path $P_1$ has no attachment
outside of the three segments incident with $v$, then
replacing $x_1Z_1v$ by $P_1$ is a T-rerouting that is
$F'$-safe, where $F'$ is as in (3),
and $P_2,x_1P_1x_2,x_2P_1v$ is a local triad.
Thus in this case the theorem holds by (3).
We may therefore assume that
there exists a path $P$ with ends $x\in V(P_1)-\{x_1,y_1\}$ and
$y\in V(S)$ such that $P$ has no internal vertex in $S\cup P_1\cup P_2$
and $y$ does not belong to any of the segments of $S$ incident with $v$.
If $y\in V(C)$, then there exists a free $S$-cross, a case we already handled,
and so we may assume not.
For $i=1,2$ let $C_i$ be the disk other than $C$ that includes $Z_i$.
Since $y_1$ belongs to the interior of $Z_2$, the only two disks it
belongs to are $C$ and $C_2$.
We may assume that $y\in V(C_2)$, for otherwise $P_1\cup P$ includes
an $S$-jump (with ends $y_1$ and $y$), in which case the theorem holds by (2).
But $C_1\cap C_2$ is equal to the third segment incident with $v$,
and hence $x_1\not\in V(C_1\cap C_2)$, and therefore $x_1\not\in V(C_2)$.
Thus $P_1\cup P$ includes an $S$-jump with ends $x_1$ and $y$,
and so the theorem holds by (2).
This \xx{proves (6)}.

Since every unstable $S$-bridge is $2$-separated from $S$ we may
apply~\refthm{planarlemma2} to the graphs $G,S$ and $H'$ to deduce
that one of the outcomes (i)--(vi) of that lemma holds.
But we may assume that (i) does not hold by (2),
we may assume that (ii) does not hold by (\xx{6}),
we may assume that (iii) does not hold, because otherwise outcome (iv)
of the present theorem holds,
we may assume that \refthm{planarlemma2}(iv) does not hold by (4),
and we may assume that \refthm{planarlemma2}(vi) does not hold,
for otherwise outcome (v) holds.
Thus we may assume that \xx{either}~\refthm{planarlemma2}(v) 
\xx{or~\refthm{planarlemma2}(vii)} holds.

\xx{%
Assume first that~\refthm{planarlemma2}(vii) holds, and let the notation 
be as in the definition of $S$-leap  as introduced
prior to~\refthm{planarlemma2}. Assume first that the degree of $v$ is three.
Let $F'$ be obtained from $F$ by removing all edges incident with $u$, where 
$\eta(u)=v$. Let $\eta'$ be obtained from $\eta$ by replacing $x_3Z_1v$ by $P_3$.
Then $\eta'$ is obtained from $\eta$ by an $F'$-safe $T$-rerouting and $x_1P_1x_2,x_2P_1y_1, P_2$
is a local triad centered at $\eta'(u)$.
Thus $F'$ and $\eta'$ satisfy the  theorem  by (3).
We may therefore assume that the degree of $v$ is at least four.
Let $\eta'$ be obtained from $\eta$ by replacing $x_1Z_1v$ by $x_1P_1x_2\cup P_2$;
then $\eta'$ is obtained by an $F'$-safe $I$-rerouting by the last axiom in the definition 
of $S$-leap.
It follows that  $F'$ and $\eta'$ satisfy the  theorem  by (1) and~(5) applied to $\eta'$ and $F'$.
}%

\xx{%
We may therefore assume that~\refthm{planarlemma2}(v) holds.}
Let the notation be as in the definition of $S$-tunnel as introduced
prior to~\refthm{planarlemma2}, and let $e_0\in E(G)$ be such that
$\eta(e_0)=W$.
Let $D$ and $D'$ be cycles in $G$ such that $\eta(D)=C$ and $\eta(D')=C'$.
If one of $B_1,B_2$ is equal to $B_{ez}$ for some $e\in E(D)-\{e_0\}$,
then such an $e$ is unique (because if $B_1\ne B_2$, then both have the
same unique attachment outside $W$),
and we denote it by $e_1$; otherwise $e_1$ is undefined.
If $e_1$ is well-defined we define $F':=F-\{e_0,e_1\}$;
otherwise we define  $F':=F-\{e_0\}$.
Let $Z':=(Z_e:e\in F')$.
Let $\eta':G\emb S'\subseteq H'$ be obtained from $\eta$ by replacing
$x_1Wy_1$ by $P_1$.
Then this rerouting is $F'$-safe, and so $Z'$ is feasible for $\eta'$
by~\refthm{reroutefeas}.
By~\refthm{unitefull} we may assume that there is a full cast for
$Z'$ and $\eta'$ in $H'$, for otherwise the theorem holds.
Let $\Gamma'=(B'_{ez}:e\in F', z\in Z_e)$ be such a cast as constructed
in the proofs of~\refthm{reroutefeas} and~\refthm{unitefull}.
Let $B'$ be the $S'$-bridge containing $P_2$ and $P_4$.
Then $B'$ is a subgraph of the union of $B_1$, $B_2$, $x_1Wy_1$ and
all $S$-bridges that have an attachment in the interior of $x_1Wy_1$.
It follows from the definition of $S$-tunnel by analyzing the
proof of~\refthm{reroutefeas} that if $B'=B'_{ez}$ for some $e\in F'$
and $z\in Z_e$, then $e\in E(D')-\{e_0\}$.
The $S'$-bridge $B'$ includes an $S'$-path $P$ with
one end say $x\in V(C)-V(W)$ and the other end say $y\in V(C')-V(W)$.
We may assume that $x\in V(\eta(e_1))$ if $e_1$ is well-defined.
In a manner similar as before, by replacing $y$ by a different vertex if
necessary, we may choose $P$ to be compatible with $\Gamma'$.
If $P$ is an $S'$-jump, then outcome (ii) holds, and so we may assume
that it is not.
Thus some disk $C''$ of $S'$ includes both $x$ and $y$.
It follows that $B'$ includes an $S'$-triad.
By~\refthm{planartriad} the triad is local; let it be centered at $v\in V(S)$.
It follows that $\eta(e_0)$ is incident with $v$, and so is $\eta(e_1)$
if $e_1$ is well-defined (by the choice of $x$).
Thus the theorem holds by (3).~\qed

We deduce the following corollary.

\newthm{apexcor}
Let $G$ be an \ifc\ triangle-free planar graph not isomorphic to the cube,
and let $F\subseteq E(G)$ be such that no two elements of $F$ belong to
the same facial cycle of $G$.
Let $H$ be a graph, and let $Z=(Z_e:e\in F)$ be a mold for $G$ in $H$.
Let $H':=H\backslash\bigcup_{e\in F} Z_e$, and let
$\eta_0:G\emb S_0\subseteq H'$ be a \he\ such that the mold $Z$ is feasible
for $\eta_0$.
If $H'$ is \ifc\ and non-planar,
then there exists a set $F'\subseteq F$ with $|F-F'|\le 1$
such that the graph $L$ determined by $G$ and $(Z_e:e\in F')$
satisfies one of the following conditions:
\item{(i)}there exist vertices $u,v\in V(L)-V(Z)$ that do not belong
to the same facial cycle of $L\backslash V(Z)$ such that
$L+uv$ is isomorphic to a minor of $H$,
\item{(ii)}there exists a facial cycle $C$ of $L-V(Z)$
and distinct vertices $u_1,u_2,v_1,v_2\in V(C)$ appearing on $C$
in the order listed such that $L+u_1v_1+u_2v_2$ is isomorphic
to a minor of $H$, and $u_iv_i \not \in E(G)$ for $i=1,2.$

\proof
Let $\eta:G\emb S\subseteq H'$, $F'$ and $Z'=(Z_e:e\in F')$ be
as in~\refthm{mainapex}.
Then $|F-F'|\le 1$, because no two edges of $F$ are cofacial.
By~\refthm{mainapex} one of (i)--(v) of that theorem holds.
But (iv) does not hold, because $H'$ is \ifc, and  (v) does not hold, because
$H'$ is not planar.
Let $\Gamma=(B_{ez}:e\in F',z\in Z_e)$ be a cast satisfying (i), (ii), or
(iii) of~\refthm{mainapex}.
If~\refthm{mainapex}(i) holds, then let $e,f\in F$ be distinct edges
such that $B_{ez}$ and $B_{fw}$ are
subgraphs of the same $S\cup Z$-bridge for some $z\in Z_e$ and $w\in Z_f$.
It follows from the proof of~\refthm{moldminor}
that $L+\hat e\hat f$ is isomorphic to a minor of $H$,
as required for (i).
If~\refthm{mainapex}(ii) holds, then (i) holds by~\refthm{compatpathminor},
and if~\refthm{mainapex}(iii) holds, then (ii) holds
by~\refthm{compatcrossminor}.~\qed

When $V(Z)$ has size one we get the following explicit version,
which is used in~[\cite{KawNorThoWolbdtw}].

\newthm{oneapex}
Let $G$ be an \ifc\ triangle-free planar graph not isomorphic to the cube,
and let $F\subseteq E(G)$ be a non-empty set
such that no two edges of $F$
are incident with the same face of $G$.
Let $G'$ be obtained from $G$ by subdividing each edge in $F$
exactly once, and let
$L$ be the graph obtained from $G'$ by adding a new vertex $v\not\in V(G')$
and joining it by an edge to all the new vertices of $G'$.
Let a subdivision of $L$ be isomorphic to a subgraph of $H$, and let
$u\in V(H)$ correspond to the vertex $v$.
If $H\backslash u$ is \ifc\ and non-planar,
then there exists an edge $e\in E(L)$ incident with $v$ such that
either
\item{(i)} there exist vertices $x,y\in V(G')$ not belonging
to the same face of $G'$ such that $(L\backslash e)+xy$ is
isomorphic to a minor of $H$, or
\item{(ii)} there exist vertices $x_1,x_2,y_1,y_2\in V(G')$ appearing
on some face of $G'$ in order such that $(L\backslash e)+x_1y_1+x_2y_2$ is
isomorphic to a minor of $H$, and $x_iy_i \not \in E(G)$ for $i=1,2$.

\proof
For $e\in F$ let $Z_e:=\{v\}$, and let $Z=(Z_e:e\in F)$.
Then $L$ is the graph determined by $G$ and $Z$.
If we identify $u$ and $v$, then $Z$ becomes a mold for $G$ in $H$.
Since a subdivision of $L$ is isomorphic to a subgraph of $H$,
the second half of~\refthm{moldminor} implies that $Z$ is feasible
for a \he\ $\eta:G\emb S\subseteq H\backslash u$.
By~\refthm{apexcor} the corollary holds.~\qed

\newsection{pinwheel} A SECOND APPLICATION

In this section we describe an application of~\refthm{apexcor}. Let $C_1$ and $C_2$ be two vertex-disjoint cycles of length $n \geq 3$ with vertex-sets $\{x_1,x_2,\ldots,x_n\}$ and $\{y_1,y_2,\ldots,y_n\}$ (in order), respectively, and let $G$ be the graph obtained from
the union of $C_1$ and $C_2$ by adding an edge joining $x_i$ and $y_i$ for each $i = 1,\ldots,n$. We
say that $G$ is \dfn{a planar ladder with $n$ rungs} and we say that $C_1$ and $C_2$ are \dfn{the rings of $G$}. Suppose now that $n=2k$ and let $W$ be a set disjoint from $V(G)$. Let $F = \{x_{2i}y_{2i} : 1 \leq i \leq k\}$. For every $e = x_{2i}y_{2i} \in F$ define $Z_e=W$. Then $Z = (Z_e : e \in F)$ is a mold for $G$ and we refer to it as a \dfn{$|Z|$-pinwheel mold}. Let $L$ be the graph determined by $G$ and $Z$. We say that $L$ is a \dfn{$|Z|$-pinwheel with $k$ vanes}.

Let $G'$ be a graph obtained from the graph $G$ described above by deleting the edges $x_1x_n$ and $y_1y_n$ and adding the edges $x_1y_n$ and $y_1x_n$. Then we say that $G'$ is \dfn{a M\"{o}bius ladder with $n$ rungs}. Let $Z$ be defined as in the previous paragraph and let $L'$  be the graph determined by $G'$ and $Z$. We say that $L'$ is a \dfn{M\"{o}bius $|Z|$-pinwheel with $k$ vanes}.

\newthm{pinwheelmold} Let $k,t$ be positive integers. Let $G$ be the planar ladder with $8(k+1)$ rungs. Let $H$ be a $(t+4)$-connected graph, and let $Z=(Z_e:e\in F)$ be a $t$-pinwheel mold for $G$ in $H$. Let $H':=H\backslash V(Z)$, and let
$\eta_0:G\emb S_0\subseteq H'$ be a \he\ such that the mold $Z$ is feasible
for $\eta_0$. Then either
\item{(i)} $H\setminus V(Z)$ is planar, or
\item{(ii)} $H$ has a minor isomorphic to a M\"{o}bius $t$-pinwheel with $k$ vanes.

\proof By \refthm{apexcor} either \refthm{pinwheelmold}(i) holds or there exists a set $F' \subseteq F$ with $|F-F'| \leq 1$ such that the graph $L$ determined by $G$ and $(Z_e:e\in F')$ satisfies one of the outcomes \refthm{apexcor}(i) and (ii). 
Let us consider the case when $L$ satisfies \refthm{apexcor}(ii), as the argument in the other case is analogous. Let $L,C,u_1,v_1,u_2,v_2$ be as in  \refthm{apexcor}(ii). If $C$ is not a ring of $G$, then it is easy to see that 
$(L\backslash V(Z))+u_1v_1+u_2v_2$ has a minor isomorphic to a M\"{o}bius ladder with $8(k+1)$ rungs, 
and by removing at most two rungs we find a 
subdivision of a M\"{o}bius $t$-pinwheel with $4k$ vanes in $L+u_1v_1+u_2v_2$.

Suppose now that $C$ is a ring of $G$. Without loss of generality we may assume that $u_1,u_2,v_1,v_2 \in \{x_i : 2(k+1) \leq i \leq 8(k+1)\}$. Then $G+u_1v_1+u_2v_2$ contains a subdivision  of a M\"{o}bius ladder with $2(k+1)$ rungs with branch vertices $x_1,x_2,\ldots,x_{2(k+1)},$ $ y_1,y_2,\ldots,y_{2(k+1)}.$ It follows that $L+u_1v_1+u_2v_2$ contains a subdivision of  a M\"{o}bius $t$-pinwheel with $k$ vanes, as desired.~\qed

Note that a M\"{o}bius $t$-pinwheel with $6t$ vanes has a minor isomorphic to $K_{t+5}$; see Figure~\reffig{mobius} for an example. Thus Lemma~\refthm{pinwheelmold} implies the following theorem, which is used in~[\cite{NorThoLinear}].

\centerline{\includegraphics[scale=1.0]{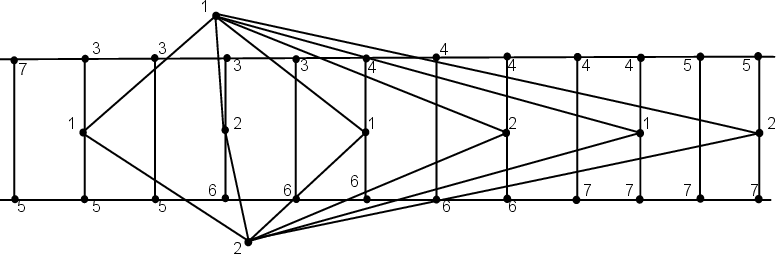}}
\bigskip
\centerline{Figure~\newfig{mobius}. $K_7$ minor in a \xx{M\"{o}bius} $2$-pinwheel with $12$ vanes.}
\bigskip

\newthm{ktpinwheel} Let $t$ be a positive integer. Let $G$ be a planar ladder with $8(6t+1)$ rungs. Let $H$ be a $(t+4)$-connected graph, and let $Z$ be a $t$-pinwheel mold for $G$ in $H$. Let $H':=H\backslash V(Z)$, and let
$\eta_0:G\emb S_0\subseteq H'$ be a \he\ such that the mold $Z$ is feasible
for $\eta_0$. Then either
\item{(i)} $H\setminus Z$ is planar, or
\item{(ii)} $H$ has a minor isomorphic to $K_{t+5}$.

\beginsection Acknowledgement

We thank Alexander Kelmans for carefully reading an earlier version
of this manuscript,  providing extensive comments and pointing out
several errors, 
\xx{and Katherine Naismith for pointing out a missing outcome in~\refthm{planarlemma2}.}
We would like to acknowledge that this paper grew out of ideas developed
in the work of the second author with Neil Robertson and Paul Seymour.

\beginsection References

\myinput{simplextbib.tex}

\medskip
\baselineskip 11pt
\vfill
\noindent
This material is based upon work supported by the National Science Foundation.
Any opinions, findings, and conclusions or recommendations expressed in
this material are those of the authors and do not necessarily reflect
the views of the National Science Foundation.
\eject

\end